\documentclass[lettersize,journal]{IEEEtran}

\usepackage{array}
\usepackage{textcomp}
\usepackage{stfloats}
\usepackage{url}
\usepackage{verbatim}
\usepackage{graphicx}
\usepackage{cite}

\usepackage{amsmath,amssymb,amsfonts}

\def\BibTeX{{\rm B\kern-.05em{\sc i\kern-.025em b}\kern-.08em
    T\kern-.1667em\lower.7ex\hbox{E}\kern-.125emX}}

\usepackage{xcolor}
\usepackage{bm}
\usepackage{algorithm}
\usepackage{algpseudocode}



\newcounter{defcounter}
\newenvironment{pequation}{%
\addtocounter{equation}{-1}
\refstepcounter{defcounter}
\renewcommand\theequation{{\bf{P\thedefcounter}}}
\begin{equation}}
{\end{equation}}

\usepackage{subfigure}

\usepackage{cases}
\usepackage{array}

\usepackage{empheq}

\usepackage{dsfont}

\newcommand{\mybold}{\text{\usefont{U}{bbold}{m}{n}1}}
\MakeRobust{\mybold}

\begin{document}

\title{Privacy-Preserving Distributed Energy Resource Control with Decentralized Cloud Computing}

\author{Xiang Huo, \IEEEmembership{Graduate Student Member, IEEE}, Mingxi Liu, \IEEEmembership{Member, IEEE}
\thanks{The authors are with the Department of Electrical and Computer Engineering, University of Utah, Salt Lake City, UT 84112 USA (e-mail: xiang.huo, mingxi.liu@utah.edu).}}





\maketitle

\begin{abstract} The rapidly growing penetration of renewable energy resources brings unprecedented challenges to power distribution networks -- management of a large population of grid-tied controllable devices encounters control scalability crises and potential end-user privacy breaches. Despite the importance, research on privacy preservation of distributed energy resource (DER) control in a fully scalable manner is lacked. To fill the gap, this paper designs a novel decentralized  privacy-preserving DER control framework that 1) achieves control scalability over DER population and heterogeneity; 2) eliminates peer-to-peer communications and secures the privacy of all participating DERs against various types of adversaries; and 3) enjoys higher computation efficiency and accuracy compared to state-of-the-art privacy-preserving methods. A strongly coupled optimization problem is formulated to control the power consumption and output of DERs, including solar photovoltaics and energy storage systems, then solved using the projected gradient method. Cloud computing and secret sharing are seamlessly integrated into the proposed decentralized computing to achieve privacy preservation. Simulation results prove the capabilities of the proposed approach in DER control applications.

\end{abstract}



\section{Introduction}

\subsection{Related Works}

\IEEEPARstart{L}{arge}-scale deployment of distributed energy resources (DERs) has proven efficacy in reducing carbon footprint and  providing grid-edge  services such as voltage control, load following, and backup power supply \cite{derancillary}. DERs, including energy storage systems (ESSs), solar photovoltaic (PV), and electric vehicles (EVs), along with other monitoring and controllable devices, can offer significant opportunities for advancing efficient, reliable, and cost-effective power grids \cite{abb_aguero2017modernizing,abb_song2012operation}. Though integrating DERs into power grids can provide multifarious benefits, such as enhanced energy efficiency and economic boost, the high penetration of DERs raises surging challenges on the scalability of existing control strategies \cite{abb_molzahn2017survey}.

To address the aforementioned challenges in large-scale DER control problems, distributed and decentralized control strategies are drawing increased attention owing to their superior scalability. For instance, a distributed  coordination method based on local droop control and consensus control was designed in \cite{abb_zeraati2016distributed} to deal with the voltage rise problem caused by the high penetration of solar PVs. Zhang \emph{et al.} in \cite{abb_zhang2021distributed} proposed an asynchronous distributed leader-follower control strategy that optimally schedules DERs to lower the voltage for peak load shaving and long-term energy saving. To reduce the communication burden, a distributed low-communication algorithm was proposed in \cite{abb_pan2021distributed} to control islanded PV-battery-hybrid systems. Though distributed methods can achieve scalability, they generically suffer from massive peer-to-peer communications. To overcome this issue, Navidi \emph{et al.} in \cite{abb_navidi2018two} developed a two-layer decentralized DER coordination architecture that can scale the solution to large networks, and no direct communication is required between local controllers. In \cite{abb_lin2017decentralized}, a decentralized stochastic control strategy was designed 
for radial distribution systems with controllable PVs and ESSs to minimize the demand balancing cost. Huo \emph{et al.}
in \cite{abb_huo2022two} proposed a decentralized shrunken primal-multi-dual subgradient algorithm with dimension reduction to achieve scalability \emph{w.r.t.} both agent population size and network dimension.

Despite the superior scalability and communication efficiency of decentralized methods, their implementation has been significantly hampered by the vulnerability to privacy breaches. Furthermore, both distributed and decentralized strategies rely heavily on mandatory communications which can disclose users' sensitive information and expose system vulnerabilities to adversaries. Differential privacy (DP) has received substantial attention in addressing privacy concerns due to its rigorous mathematical formulation \cite{abb_dwork2006calibrating}. DP-based methods add persistent randomized perturbations to the datasets, constraints, or objective functions for privacy preservation. In \cite{abb_dong2018privacy}, a DP-based aggregation algorithm is proposed to compensate for solar power fluctuations and protect users' personal information. Han \emph{et al.} in \cite{abb_han2016differentially} developed a distributed optimization algorithm based on DP to preserve the privacy of the participating agents. 
Gough \emph{et al.} in \cite{abb_gough2021preserving} designed an innovative DP-compliant algorithm to ensure that the data from consumers' smart meters are protected. Despite the success in privacy preservation, DP-based methods inevitably suffer from accuracy loss due to the added perturbations.


In contrast, encryption-based strategies achieve privacy preservation with high accuracy by encrypting the original data into cyphertexts, and only those holding private keys can decrypt the cyphertexts. Lu \emph{et al.} in \cite{abb_lu2012eppa} proposed an efficient and privacy-preserving aggregation scheme for smart grid communications, in which the data is encrypted by Paillier cryptosystem. In \cite{abb_mohammadali2021privacy}, a privacy-preserving and fault-tolerant scheme was designed based on homomorphic cryptosystem to achieve secure aggregation of metering data. Similarly, Cheng \emph{et al.} in \cite{abb_cheng2021homomorphic} proposed a novel private collaborative distributed energy management system based on homomorphic encryption to solve the privacy issues in distribution systems and microgirds. Despite the high accuracy, the drawback of encryption-based methods lies in the prevalent computing overhead caused by encryption and decryption. Other hardware-integrated privacy-preserving methods, e.g., garbled circuit \cite{abb_wang2018privacy,abb_gilad2019secure}, are deficient in flexibility and  uneconomic due to the hardware cost.

Secret sharing (SS) \cite{abb_shamir1979share} is a lightweight cryptographic method that can securely distribute a secret among a group of participants. Each participant will be allocated a share of the secret, and only through the collaboration of certain participants where the number of participants is greater than a threshold can the secret be reconstructed from their shares. Adopting SS, Nabil \emph{et al.} in \cite{nabil2019ppetd} designed an SS-based detection scheme to identify malicious consumers who steal electricity, in which system operators only collect masked meter readings from the consumers to avoid privacy violation. 
In \cite{abb_huo2022distributed}, an SS-based EV charging control protocol was developed to achieve privacy-preserving EV charging control for overnight valley filling. 
Compared with encryption-based strategies, SS-based methods can preserve privacy while avoiding the heavy computational load. Despite the superiority, few research studied the integration of SS into DER control due to the highly complex distribution network structure, large DER population, and lack of theoretical support in privacy guarantees. To fill these gaps, this paper designs a novel SS-based privacy-preserving algorithm that merits high efficiency, security, and accuracy for large-scale DER control problems.

\vspace{-4mm}
\subsection{Statement of Contributions}



The contribution of this paper is three-fold: 1) We propose a novel decentralized privacy-preserving algorithm that concurrently achieves scalability and privacy in large-scale DER control. To the best of our knowledge, this is the first paper that proposes a decentralized SS-based algorithm for DER privacy preservation, in which decentralized solutions, privacy guarantees, and rigorous security proofs are provided; 2) The proposed method eliminates the frequent peer-to-peer communications and secures the privacy of the participating DERs against various types of adversaries. The designed framework serves as a benchmark for secure and scalable DER control. 3) Compared to state-of-the-art approaches, the proposed method can achieve lower computational overhead and identically accurate solutions as the non-privacy-concerned algorithms. 

The rest of this paper is organized as follows: In Section \ref{Problem_formulation}, we construct the models of distribution networks, PVs, and ESSs, then formulate the DER control problem into a constrained optimization problem. Section \ref{Decentralized_optimization} derives the decentralized solution via the projected gradient method and presents the corresponding DER aggregation and control strategies. The SS-based privacy-preserving DER control algorithm and privacy analyses are provided in Section \ref{Privacy_preserving}. We give simulation results and analyses in Section \ref{Simulation_results}. Section \ref{Conclusion} concludes this paper.

\section{Problem Formulation}
\label{Problem_formulation}

\subsection{Branch Flow Model} 
\begin{figure}[!htb]
\vspace{-2mm}
    \centering
 \includegraphics[width=0.41\textwidth, trim={0.1cm 0cm 0.0cm 0.1cm},clip]{ 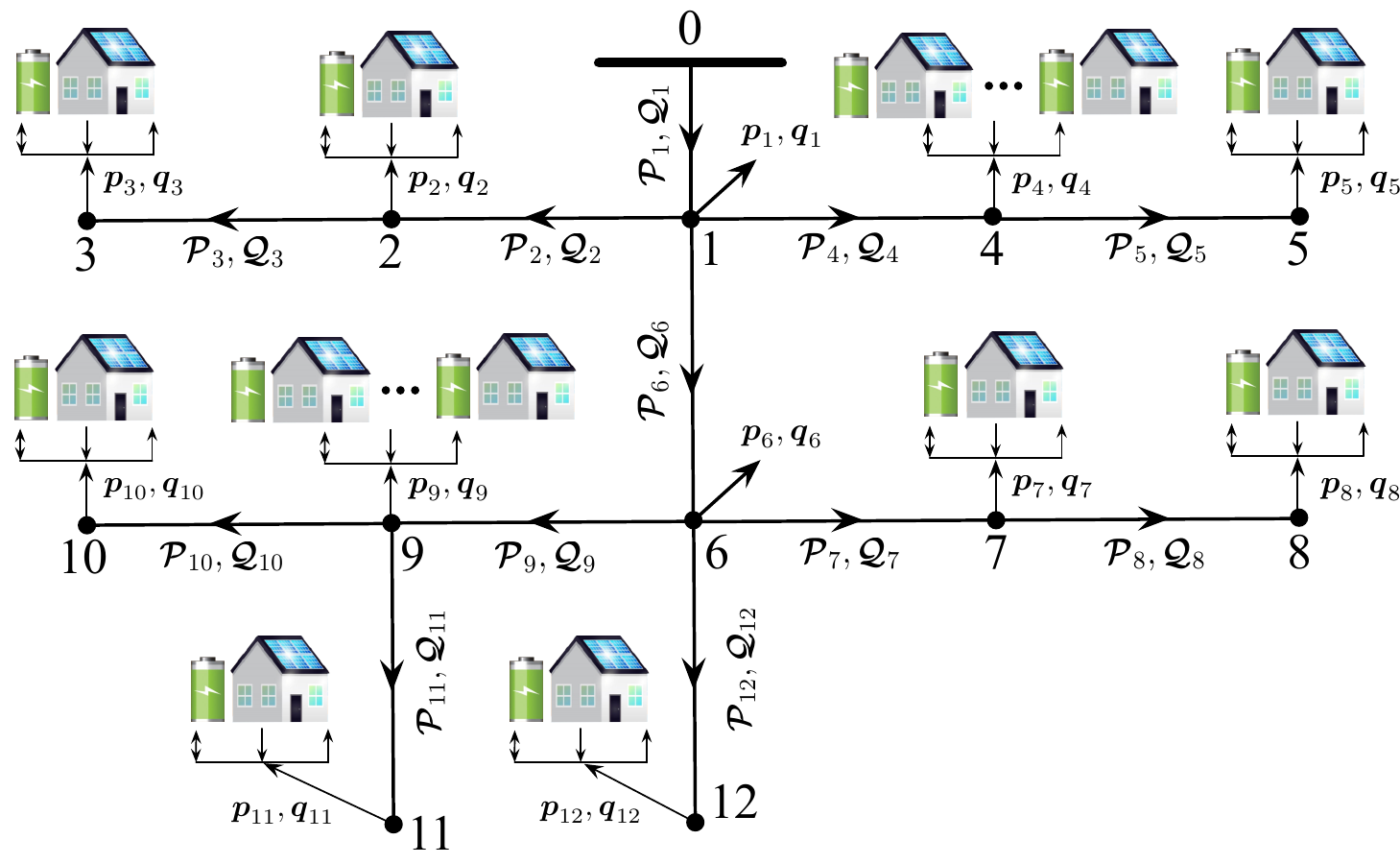}
 \vspace{-2mm}
    \caption{A radial 13-bus distribution network connected with rooftop solar PVs and ESSs.}
    \label{13bus}
     \vspace{-2mm}
\end{figure}

Consider an $n$-bus radial distribution network where $\mathbb{B} = \{0,1,\ldots,n \}$ denotes the set of buses. Let $l_{ij}$ denote the line segment connecting buses $i$ and $j$, $\mathbb{L} = \{1,\ldots, h \}$ denote the set of lines, $\mathbb{C}_j$ denote the set of bus $j$'s child buses, $V_j$ denote the voltage magnitude at bus $j$, $\mathcal{P}_{ij}$ and $\mathcal{Q}_{ij}$ denote the active and reactive power flow from bus $i$ to bus $j$, respectively, and $r_{ij}$ and $x_{ij}$ be the resistance and reactance of line $l_{ij}$, respectively. For bus $j$, let  $p_j^c$ and $q_j^c$ denote the active and reactive power consumptions, respectively, and $p_j^g$ and $q_j^g$ denote its active and reactive power generations, respectively. 
To simplify the network model, a nonlinear DistFlow model \cite{abb_baran1989optimal} can be linearized to the LinDistFlow model by omitting the higher order terms with negligible 
error \cite{abb_farivar2013equilibrium}. Therefore, this paper adopts the LinDistFlow model, represented as
\begin{subequations}
\begin{align}
    \mathcal{P}_{ij} - \sum_{u \in \mathbb{C}_j} \mathcal{P}_{ju} &= p_j^c  - p_j^g \label{active_flow}\\
    \mathcal{Q}_{ij} - \sum_{u \in \mathbb{C}_j} \mathcal{Q}_{ju} &= q_j^c  - q_j^g  \\
    V_i^2 - V_j^2 &= 2(r_{ij} \mathcal{P}_{ij} + x_{ij}\mathcal{Q}_{ij}).
\end{align}
\end{subequations}
A radial 13-bus distribution network connected with rooftop solar PVs and ESSs is shown in Fig. \ref{13bus} and will be used as an example throughout this paper.

In this paper, one control objective is to minimize the total power loss of the distribution network by controlling the dynamics of PVs and ESSs, which is approximated by
\begin{equation}
f_1(\bm{p}_1^g,\ldots,\bm{p}_n^g) = \sum_{l_{ij} \in \mathbb{L}} r_{ij} \left
(\frac{\| \bm{\mathcal{P}}_{ij} \|^2_2 + \| \bm{\mathcal{Q}}_{ij} \|^2_2}{V_0^2}\right)\vspace{-2mm}
\end{equation}
where $V_0$ denotes the nominal voltage magnitude,  $\bm{p}_j^g$, $\bm{\mathcal{P}}_{ij}$, and 
$\bm{\mathcal{Q}}_{ij} \in \mathbb{R}^{T}$ are augmented vectors of $p_j^g$, $\mathcal{P}_{ij}$, and 
$\mathcal{Q}_{ij}$ across $T$ time intervals, respectively. Note that we only consider active power loss and assume reactive power flows $\bm{\mathcal{Q}}_{ij}$ to be constant vectors. Though the reactive power loss is not included here for simplicity, it can be added without affecting algorithm design. The active power flows are constrained by 
\begin{equation}
  \bm{0} \leq \bm{\mathcal{P}}_{ij} \leq \overline{\bm{\mathcal{P}}}_{ij}
  \label{flow_limit}
\end{equation} where $\overline{\bm{\mathcal{P}}}_{ij}$ denotes the maximum active power flow limit.

\vspace{-4mm}
\subsection{Solar Photovoltaic} Let $\mathbb{V}$ denote the set of in total $\mathcal{V}$ solar PVs. During $T$ time intervals of a day, the active power injection $\tilde{\bm{p}}_{\nu} \in \mathbb{R}^{T}$  from the $\nu$th PV inverter should satisfy 
\begin{equation}
\bm{0} \leq \tilde{\bm{p}}_{\nu} \leq \bm{\overline{p}}^v_{\nu}
\label{pv_limit}
\end{equation}
where $\bm{\overline{p}}^v_{\nu}$ denotes the maximum active power injection and is assumed to be known by the forecast. Herein, the curtailment cost can be calculated  by \cite{abb_li2019distributed} 
\begin{equation}
f_2(\tilde{\bm{p}}_{\nu}) = \| \tilde{\bm{p}}_{\nu} - \bm{\overline{p}}^v_{\nu}\|^2_2.
    \label{6ss}
\end{equation}

\vspace{-4mm}
\subsection{Energy Storage System}
Let $\mathbb{S}$ denote the set of $\mathcal{E}$ ESSs. The charging/discharging power $\hat{\bm{p}}_{\sigma} \in \mathbb{R}^{T}$  of the $\sigma$th ESS is constrained by 
\begin{equation}
-\underline{\bm{p}}_\sigma^{s} \leq \hat{\bm{p}}_{\sigma} \leq \overline{\bm{p}}_\sigma^{s}
\label{charge_limit}
\end{equation}
where $\underline{\bm{p}}_\sigma^{s}$ and $\overline{\bm{p}}_\sigma^{s}$ denote the maximum discharging and charging power, respectively. Let $s^0_\sigma$ denote the initial state of charge (SoC) of the $\sigma$th ESS and $\bm{H}_\sigma \triangleq [s^0_\sigma,\ldots,s^0_\sigma]^{\mathsf{T}}\in \mathbb{R}^{T}$.  Aggregate the charging/discharging power across $T$ time intervals, then the capacity of the $\sigma$th ESS is constrained by 
\begin{equation}
 \underline{\bm{p}}_\sigma^{a} \leq \bm{H}_\sigma +  \bm{A} \hat{\bm{p}}_{\sigma} \Delta T \leq
 \overline{\bm{p}}_\sigma^{a} \label{state_limit}
\end{equation}
where $\underline{\bm{p}}_\sigma^{a}$ and $\overline{\bm{p}}_\sigma^{a}$ denote its lower and upper capacity bounds, respectively, $\Delta T$ denotes the sampling time, and the aggregation matrix $\bm{A}$ is lower triangular consisting of ones and zeros, i.e., element $\bm{A}_{\hat{\imath},\hat{\jmath}} = 1 ~ \text{if}~ \hat{\imath} \geq \hat{\jmath}$, element $\bm{A}_{\hat{\imath},\hat{\jmath}} = 0 ~ \text{if}~ \hat{\imath} < \hat{\jmath}, \forall \hat{\imath},\hat{\jmath} = 1,\ldots, T$. Therefore, the SoCs of ESS $\sigma$ during $T$ time slots are obtained by aggregating the charging/discharging power using $\bm{A}$.


Furthermore, the $\sigma$th ESS's degradation cost is calculated in terms of the smoothness of charging and discharging by \cite{abb_forman2013optimization}
\begin{equation}
f_3(\hat{\bm{p}}_{\sigma}) = \| \bm{B}\hat{\bm{p}}_{\sigma} \|^2_2.
    \label{3s}
\end{equation}
where $\bm{B}$ calculates discharging/charging differences between adjacent times, i.e., $\bm{B}_{\hat{\imath},\hat{\imath}} = 1$, $\forall \hat{\imath} = 1,\ldots, T$, $\bm{B}_{\hat{\imath},\hat{\imath}+1} = -1, \forall \hat{\imath} = 1,\ldots, T-1$, and all other elements are zeros.

 
\vspace{-2mm}
\subsection{Problem Formulation}

The optimization problem is then formulated to minimize the summation of total active power loss, PV curtailment cost, and ESS degradation cost within the distribution network as
\begin{pequation} \label{p1}
\begin{aligned}
& \underset{\tilde{\bm{p}},\, \hat{\bm{p}}}{\text{min}} & &  \delta_1 f_1(\bm{p}^g) + \sum_{\nu=1}^{\mathcal{V}} \delta_2f_2( \tilde{\bm{p}}_{\nu}) + \sum_{\sigma=1}^{\mathcal{E}}\delta_3f_3(\hat{\bm{p}}_{\sigma}) \\
& \: \text{s.t.} & & \eqref{active_flow}, \eqref{flow_limit}, \eqref{pv_limit}, \eqref{charge_limit}, \eqref{state_limit}
\end{aligned}
\end{pequation}where $\tilde{\bm{p}} = [\tilde{\bm{p}}_1^{\mathsf{T}}, \ldots, \tilde{\bm{p}}_{n}^{\mathsf{T}}]^{\mathsf{T}}$, $\hat{\bm{p}} = [\hat{\bm{p}}_1^{\mathsf{T}}, \ldots, \hat{\bm{p}}_{n}^{\mathsf{T}}]^{\mathsf{T}}$, $\bm{p}^g = [{\bm{p}_{1}^g}^{\mathsf{T}},\ldots, {\bm{p}_{n}^g}^{\mathsf{T}}]^{\mathsf{T}}$, and $\delta_\alpha$ denotes the cost coefficient associated with the objective function $f_\alpha(\cdot)$. Note that the cost coefficients are constants that allow flexible adjustments on the weights of the global and local objective functions and regulate different units.


\section{Decentralized Optimization}
\label{Decentralized_optimization}
\subsection{Projected Gradient Method}

This paper achieves scalability in solving \eqref{p1} via projected gradient method (PGM). PGM decomposes a centralized optimization problem into local optimizations at agents, resulting in a paralleled computing structure. Let $\mathbb{M} = \{1,\ldots,m\}$ denote the set of agents, e.g., buses or DERs, who work  cooperatively in solving \eqref{p1}. In this setting, the $\kappa$th agent updates its decision variable $\bm{x}_\kappa$ using PGM by 
\begin{equation} \label{eq1}
    \bm{x}_\kappa^{(\ell+1)} = \mathbb{P}_{\mathbb{X}_\kappa}[\bm{x}_\kappa^{(\ell)} - \gamma_\kappa^{(\ell)} \Phi_\kappa (\bm{x}^{(\ell)})]
\end{equation} where $\ell$ denotes the iteration number, $\bm{x}^{(\ell)} = [{\bm{x}_1^{(\ell)}}^{\mathsf{T}},\ldots,{\bm{x}_m^{(\ell)}}^{\mathsf{T}}]^{\mathsf{T}}$ includes all decision variables, i.e., $\tilde{\bm{p}}_{\nu}$ and $\hat{\bm{p}}_{\sigma}$ in problem \eqref{p1}, $\gamma_k^{(\ell)}$ denotes the step size, $\Phi_\kappa(\cdot)$ denotes the gradient of the Lagrangian \emph{w.r.t.} $\bm{x}_\kappa^{(\ell)}$, and $\mathbb{P}_{\mathbb{X}_\kappa}[\cdot]$ denotes the projection operation onto set $\mathbb{X}_\kappa$.

In \eqref{p1}, the local constraint of the $\nu$th PV in \eqref{pv_limit} and local constraints of the $\sigma$th ESS in \eqref{charge_limit} and  \eqref{state_limit} can be represented by two feasible sets $\mathbb{P}_\nu^v$ and $\mathbb{P}_\sigma^e$ as 
\begin{subequations}
\begin{align}
\mathbb{P}_\nu^v &\triangleq \{ \tilde{\bm{p}}_{\nu} | \: \bm{0} \leq \tilde{\bm{p}}_{\nu} \leq \bm{\overline{p}}^v_{\nu} \}\\
\mathbb{P}_\sigma^e &\triangleq
\{ \hat{\bm{p}}_{\sigma} | -\underline{\bm{p}}_\sigma^{s} {\leq} \hat{\bm{p}}_{\sigma} {\leq} \overline{\bm{p}}_\sigma^{s}{,}~ \:  \underline{\bm{p}}_\sigma^{a} {\leq} \bm{H} {+}  \bm{A} \hat{\bm{p}}_{\sigma} \Delta T \leq \overline{\bm{p}}_\sigma^{a}\}.
\end{align}
\end{subequations}

In what follows, aiming at reducing the number of coupling terms, we rewrite the networked constraints in \eqref{active_flow} and \eqref{flow_limit} to a single inequality constraint based on the network topology. To this end, we first represent the active power flows in \eqref{active_flow} through active power generations of each bus using
\begin{equation}
\bm{p}_i = \tilde{\bm{p}}_i - \hat{\bm{p}}_i - \bm{p}_i^c \label{total}
\end{equation}
where $\bm{p}_i$ denotes the aggregated active power generation at bus $i$,  $\tilde{\bm{p}}_i = \sum_{\nu=1}^{\mathcal{V}_i} \tilde{\bm{p}}_{\nu}$ and $\hat{\bm{p}}_i = \sum_{\sigma=1}^{\mathcal{E}_i}\hat{\bm{p}}_{\sigma}$ denote the aggregated active power of all PVs and ESSs that are connected at bus $i$, respectively.  $\mathcal{V}_i$ and $\mathcal{E}_i$ denote the total number of PVs and ESSs connected at bus $i$, respectively.


For the $\iota$th line flow $\bm{\mathcal{P}}_\iota$ in the distribution network, the from-bus is defined by the bus where the flow begins, and the to-bus set is defined by the set of buses that the $\iota$th line flow travels to till reaching the edge of the distribution network. 
Let $\bm{Z} \in \mathbb{R}^{n \times n}$ denote the adjacency matrix of the distribution network and $\bm{Z}_\iota$ denote the $\iota$th row of $\bm{Z}$ that represents the adjacency vector of the $\iota$th line flow. Let $\bm{Z}_{\iota}(i)$ denote the $i$th element of $\bm{Z}_{\iota}$, and $\bm{Z}_{\iota}(i) = 1$ if the $\iota$th power flow has bus $i$ as a to-bus, e.g., $\bm{Z}_9 = [0,0,0,0,0,0,0,0,1,1,1,0]$. Then, the power flows in the distribution network can be represented by $\bm{Z}$. Expand $\bm{Z}$ across $T$ time slots, we have 
\begin{equation}
     \Tilde{\bm{Z}} = \left[\arraycolsep=1.4pt\def\arraystretch{1.5}
    \begin{array}{cccc} 
\bm{Z}_{1}(1) \bm{I} &  \bm{Z}_{1}(2) \bm{I} & \cdots &\bm{Z}_{1}(n) \bm{I}\\ 
\vdots & \vdots &  & \vdots \\
\bm{Z}_{n}(1) \bm{I} &  \bm{Z}_{n}(2) \bm{I} & \cdots &\bm{Z}_{n}(n) \bm{I}
    \end{array}\right] 
\end{equation} where $\bm{I} {\in} \mathbb{R}^{T\times T}$ denotes the identity matrix and $\Tilde{\bm{Z}} \in \mathbb{R}^{nT \times nT}$. 

In what follows, let $\Tilde{\bm{P}} \in \mathbb{R}^{nT}$ denote the aggregated active power generations defined in \eqref{total} from all buses, we have
\begin{equation}
    \Tilde{\bm{P}} {=} \left[\begin{array}{c} 
    \bm{p}_1\\ 
    \vdots \\ 
    \bm{p}_n
    \end{array}\right] {=} \left[\arraycolsep=1.4pt\def\arraystretch{1.5}\begin{array}{c} 
    \sum_{\nu=1}^{\mathcal{V}_1} \tilde{\bm{p}}_{\nu} - \sum_{\sigma=1}^{\mathcal{E}_1}\hat{\bm{p}}_{\sigma} - \bm{p}_1^c\\ 
    \vdots \\ 
    \sum_{\nu=\mathcal{V}_{n-1} + 1}^{\mathcal{V}_n}\tilde{\bm{p}}_{\nu} - \sum_{\sigma=\mathcal{E}_{n-1} + 1}^{\mathcal{E}_n}\hat{\bm{p}}_{\sigma} - \bm{p}_n^{c}
    \end{array}
    \right]. 
\end{equation}
Furthermore, $\Tilde{\bm{P}}$ can be rewritten compactly as 
\begin{equation}
    \Tilde{\bm{P}} = \sum_{i=1}^{n}\bm{\Delta}_i\left(\tilde{\bm{p}}_i - \hat{\bm{p}}_i - \bm{p}_i^c\right) \label{16s}
\end{equation}
where $\bm{\Delta}_i$ denotes the aggregation matrix whose $i$th block is represented by the identity matrix $\bm{I}$, and all other blocks are zeros, e.g., $\bm{\Delta}_1 = [
    \bm{I}, 
     \bm{0}, 
    \ldots, 
    \bm{0}]^\mathsf{T} \in \mathbb{R}^{nT \times T}.$ Then, the active power flow of the $\iota$th line can be calculated by 
 \begin{equation}
   \bm{\mathcal{P}}_{\iota} = \Tilde{\bm{Z}}_{\iota}\Tilde{\bm{P}}. \label{global_constraint}
 \end{equation} 
Consequently, the power flow limit constraint in \eqref{flow_limit} becomes
\begin{equation}
  \bm{0} \leq \Tilde{\bm{Z}}_{\iota}\Tilde{\bm{P}} \leq \overline{\bm{\mathcal{P}}}_{\iota}.
  \label{17sss}
\end{equation}
Therefore, problem \eqref{p1} can be  written into
\begin{pequation} \label{p2}
\begin{aligned}
& \underset{\tilde{\bm{p}},\, \hat{\bm{p}}}{\text{min}} & &  \delta_1 f_1(\bm{p}^g) + \sum_{\nu=1}^{\mathcal{V}} \delta_2f_2( \tilde{\bm{p}}_{\nu}) + \sum_{\sigma=1}^{\mathcal{E}}\delta_3f_3(\hat{\bm{p}}_{\sigma}) \\
& \: \text{s.t.} & & \bm{p}_{\nu} \in \mathbb{P}_\nu^v, \ \forall \nu \in \mathbb{V} \\
& \: & & \bm{p}_{\sigma} \in \mathbb{P}_\sigma^e, \ \forall \sigma \in \mathbb{S}\\
& \: & &    \bm{0} \leq \Tilde{\bm{Z}}_{\iota}\Tilde{\bm{P}} \leq \overline{\bm{\mathcal{P}}}_{\iota}, \forall \iota \in \mathbb{L} 
\end{aligned}
\end{pequation}

The optimization problem in  \eqref{p2} seeks to find the optimal decision variables, i.e., charging and discharging power $\tilde{\bm{p}}_{\sigma}$'s of the ESSs and the active power injection $\hat{\bm{p}}_{\nu}$'s of the PVs. In what follows, we focus on solving \eqref{p2} through a decentralized fashion based on  PGM defined in \eqref{eq1}. To solve \eqref{p2} via PGM, we firstly derive its relaxed Lagrangian as \begin{align}
    \mathcal{L}(\tilde{\bm{p}}, \hat{\bm{p}}, \bm{\mu}_l, \bm{\mu}_u) &= \delta_1 f_1(\bm{p}^g) + \sum_{\nu=1}^{\mathcal{V}} \delta_2f_2( \tilde{\bm{p}}_{\nu}) + \sum_{\sigma=1}^{\mathcal{E}}\delta_3f_3(\hat{\bm{p}}_{\sigma})\nonumber\\
     & ~~~ {+} \sum_{\iota=1}^{L}\bm{\mu}_{u\iota}^{\mathsf{T}} (\Tilde{\bm{Z}}_{\iota}\Tilde{\bm{P}}-\overline{\bm{\mathcal{P}}}_{\iota})
    {-} \sum_{\iota=1}^{L}\bm{\mu}_{l\iota}^{\mathsf{T}} \Tilde{\bm{Z}}_{\iota}\Tilde{\bm{P}}
\label{Lagrangian}
\end{align}
where $\bm{\mu}_l = [\bm{\mu}_{l1}^{\mathsf{T}},\ldots,\bm{\mu}_{lL}^{\mathsf{T}}]^{\mathsf{T}}$ and $\bm{\mu}_u = [\bm{\mu}_{u1}^{\mathsf{T}},\ldots,\bm{\mu}_{uL}^{\mathsf{T}}]^{\mathsf{T}}$, $\bm{\mu}_{l\iota}$ and $\bm{\mu}_{u\iota}$ denote the dual variables associated with lower and upper power flow limits of the line $\iota$, respectively.


Suppose $\tilde{\bm{p}}_{\nu}$ and $\hat{\bm{p}}_{\sigma}$ are decision variables of the $\nu$th PV and $\sigma$th ESS  connected at bus $i$, respectively. Take the subgradients of \eqref{Lagrangian} \emph{w.r.t.} the primal variables $\tilde{\bm{p}}_{\nu}$ and $\hat{\bm{p}}_{\sigma}$, we have  
\begin{subequations} \label{22}
\begin{align}
    \nabla_{\tilde{\bm{p}}_{\nu}} \mathcal{L}(\cdot) &= 2 \delta_2 (\tilde{\bm{p}}_{\nu} - \bm{\overline{p}}_\nu^{v})
    + \frac{2\delta_1}{V_0^2}\sum_{\iota=1}^{L} r_\iota (\Tilde{\bm{Z}}_{\iota} \bm{\Delta}_i)^{\mathsf{T}} (\Tilde{\bm{Z}}_{\iota}\Tilde{\bm{P}}) \nonumber\\
    &~~~ + \sum_{\iota=1}^{L} ( \Tilde{\bm{Z}}_{\iota}\bm{\Delta}_i)^{\mathsf{T}} (\bm{\mu}_{u\iota} - \bm{\mu}_{l\iota})\label{22a}\\
\nabla_{\hat{\bm{p}}_{\sigma}} \mathcal{L}(\cdot) &= 2 \delta_3 \hat{\bm{p}}_{\sigma}
    - \frac{2\delta_1}{V_0^2}\sum_{\iota=1}^{L} r_\iota (\Tilde{\bm{Z}}_{\iota} \bm{\Delta}_i)^{\mathsf{T}} (\Tilde{\bm{Z}}_{\iota}\Tilde{\bm{P}}) \nonumber\\
    &~~~ - \sum_{\iota=1}^{L} ( \Tilde{\bm{Z}}_{\iota}\bm{\Delta}_i)^{\mathsf{T}} (\bm{\mu}_{u\iota} - \bm{\mu}_{l\iota}). \label{22b}
\end{align}
\end{subequations}

Without affecting the efficacy of the algorithm design, we assume all power lines have the same resistance $\bar{r}$ for the simplicity of presentation, herein \eqref{22} becomes
\begin{subequations} \label{22ss}
\begin{align}
    \nabla_{\tilde{\bm{p}}_{\nu}} \mathcal{L}(\cdot) &= 2 \delta_2 (\tilde{\bm{p}}_{\nu} - \overline{\bm{p}}_\nu^{v})
    + \bar{\delta}_1 \bm{\pi}_{i} \Tilde{\bm{P}} + \bm{\psi}_i (\bm{\mu}_{u} - \bm{\mu}_{l})\label{22as}\\
\nabla_{\hat{\bm{p}}_{\sigma}} \mathcal{L}(\cdot) &= 2 \delta_3 \hat{\bm{p}}_{\sigma}
    - \bar{\delta}_1 \bm{\pi}_{i}\Tilde{\bm{P}} - \bm{\psi}_{i} (\bm{\mu}_{u} - \bm{\mu}_{l}) \label{22bs}
\end{align}
\end{subequations} where $\bar{\delta}_1 = \frac{2\delta_1}{V_0^2} \bar{r}$, $\bm{\pi}_i = \sum_{\iota=1}^{L} (\Tilde{\bm{Z}}_{\iota} \bm{\Delta}_i)^{\mathsf{T}} \Tilde{\bm{Z}}_{\iota}$, and $\bm{\psi}_i$ denotes the $i$th column block of $\Tilde{\bm{Z}}$.


The detailed derivation of the Lagrangian subgradients in \eqref{22ss} can be found in \textsc{Appendix} \ref{Derivation_proof}. 

Therefore, based on the calculated subgradients in \eqref{22}, at the $\ell$th iteration, the $\nu$th PV and the $\sigma$th ESS can update their decision variables using PGM by 
\begin{subequations}\label{24}
\begin{align}
\tilde{\bm{p}}_{\nu}^{(\ell+1)} & =
\Pi_{\mathbb{P}_\nu^v}\left(\tilde{\bm{p}}_{\nu}^{(\ell)}-\alpha_{\nu,\ell}^{v} \nabla_{\tilde{\bm{p}}_{\nu}}  \mathcal{L}^{(\ell)}\left(\cdot \right)\right) \label{eq_24a}\\
\hat{\bm{p}}_{\sigma}^{(\ell+1)} & =
\Pi_{\mathbb{P}_\sigma^e}\left(\hat{\bm{p}}_{\sigma}^{(\ell)}-\alpha_{\sigma,\ell}^{e} \nabla_{\hat{\bm{p}}_{\sigma}}    \mathcal{L}^{(\ell)}\left(\cdot \right)\right) \label{eq_24b}
\end{align}
\end{subequations}
where $\alpha_{\nu,\ell}^{v}$ and $\alpha_{\sigma,\ell}^{e}$ denote the primal step sizes of the $\nu$th PV and the $\sigma$th ESS, respectively, $\mathcal{L}^{(\ell)}\left(\cdot \right)$ denotes the calculated Lagrangian in \eqref{Lagrangian} at the $\ell$th iteration. The dual variables can be updated similarly using PGM.

\subsection{DER Aggregation and Control}

In PGM iterations, the $i$th agent needs to calculate $\Phi_i (\bm{x}^{\ell})$ in \eqref{eq1} where the decision variables $\bm{x}_i$'s from all other agents are required. As indicated in \eqref{22ss},  calculating subgradients $\nabla_{\tilde{\bm{p}}_{\nu}}\mathcal{L}(\cdot)$ and $\nabla_{\hat{\bm{p}}_{\sigma}}\mathcal{L}(\cdot)$ indeed requires the decision variables $\Tilde{\bm{P}}$ from all the agents. Specifically, the calculation of subgradients in  \eqref{22as} and \eqref{22bs} are coupled through
\begin{equation}
  \mathcal{C} =  \mathcal{C}_p +  \mathcal{C}_d  =\bar{\delta}_1 \bm{\pi}_i \Tilde{\bm{P}} + \bm{\psi}_i (\bm{\mu}_{u} - \bm{\mu}_{l}) \label{27sss}
\end{equation}
where $\mathcal{C}_p$ and $\mathcal{C}_d$ denote the coupling terms associated with the primal and dual variables, respectively.  

To clearly demonstrate the information exchange needs in subgradient calculation, we exemplify the primal update of the $\hat{\nu}$th PV connected at bus 2. The $\hat{\nu}$th PV can update its decision variable $\tilde{\bm{p}}_{\hat{\nu}}$ using the subgradient in \eqref{22as} which is
\begin{equation}
\nabla_{\tilde{\bm{p}}_{\hat{\nu}}} \mathcal{L}(\cdot) = 2 \delta_2 (\tilde{\bm{p}}_{\hat{\nu}} - \bm{\overline{p}}_{\hat{\nu}}^{v})
    + \sum_{\iota=1}^{2} \left(\bar{\delta}_1\bm{\pi}_{\iota}\Tilde{\bm{P}} + \bm{\mu}_{u\iota} + \bm{\mu}_{l\iota} \right)\label{26ss}
\end{equation}
where $\bm{\pi}_{1}\Tilde{\bm{P}} = \sum_{i=1}^{n}\bm{p}_i$ and $\bm{\pi}_{2}\Tilde{\bm{P}} = \bm{p}_2 + \bm{p}_3$. Therefore, the $\hat{\nu}$th PV requires the active power generations $\bm{p}_i, \forall i=1,\ldots,n$ from all buses to conduct the update in \eqref{eq_24a}.

\begin{figure}[!htb]
\vspace{-2mm}
    \centering
\includegraphics[width=0.48\textwidth, trim={0cm 0cm 0cm 0cm},clip]{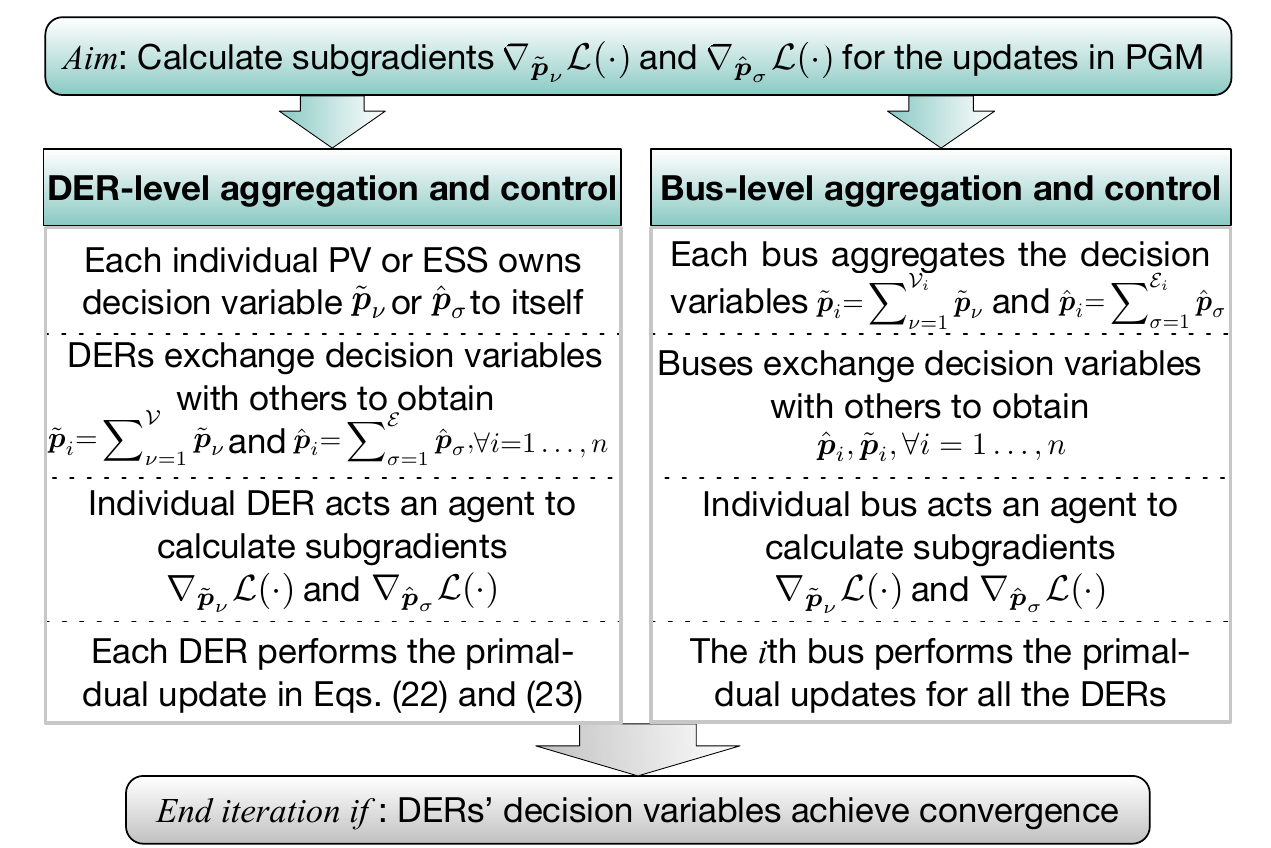}
\vspace{-2mm}
  \caption{Aggregation and control of DERs via bus-level and DER-level architectures.}
\label{bus_der_level_control}
 \vspace{-1mm}
\end{figure}

Based on the above observations, two different aggregation and control strategies, i.e.,   \textit{Bus-level aggregation and control} and \textit{DER-level aggregation and control}, can be applied as shown in Fig. \ref{bus_der_level_control}.  In bus-level aggregation and control, the $i$th bus (agent) aggregates the decision variables $\tilde{\bm{p}}_i =    \sum_{\nu=1}^{\mathcal{V}_i}\tilde{\bm{p}}_{\nu}$ and $\hat{\bm{p}}_i = \sum_{\sigma=1}^{\mathcal{E}_i}\hat{\bm{p}}_{\sigma}$ where only aggregated decision variables are transmitted and used for the primal updates. In contrast, DER-level control strategies require each DER to act as an agent and receive all data of others that is demanded for updates in \eqref{24}. However, due to the large number of DERs connected to the distribution network, DER-level control can suffer from massive data exchange and heavy local computation. Therefore, we adopt the bus-level aggregation and control scheme which is more computing and communicating efficient. We will later show that the proposed privacy-preserving algorithm can be readily extended to the DER-level control (See \textbf{Remark 1} for details). 

Apart from scalability and efficiency, the inevitable private information exposure in both bus-level and DER-level methods raises fundamental privacy concerns, e.g., the electrical load can reveal sensitive business activities and/or customer's daily routines. 
To address the privacy concerns, we will develop a novel SS-based algorithm to achieve secure information exchange in executing \eqref{24}.
 

\section{SS-Based Privacy-Preserving DER Control}
\label{Privacy_preserving}
\subsection{Real Number to Integer Quantization}

Note that the SS scheme requires modular arithmetic instead of real arithmetic. However, decentralized optimization genetically requires real number calculations, e.g., real decision variables and parameters. Therefore, a real number to integer transformation is needed to integrate SS into decentralized optimization. We adopt the fixed-point number quantization  \cite{abb_daru2019encrypted} to map the real numbers onto the integer space and the fixed-point real-number set is defined by 
\begin{equation}
    \mathbb{Q}_{\theta, \gamma, \zeta } {\triangleq} \left\{-\theta^{\gamma},-\theta^{\gamma}+\theta^{-\zeta }, \ldots, \theta^{\gamma}-2 \theta^{-\zeta }, \theta^{\gamma}-\theta^{-\zeta }\right\}
    \label{27ss}
\end{equation}
where $\theta \in \mathbb{N}_{1+}$ denotes the basis, $\gamma \in \mathbb{N}$ denotes the magnitude, and $\zeta  \in \mathbb{N}$ denotes the resolution. Therefore, by defining a 
surjective mapping $m(\cdot): \mathbb{R} \mapsto \mathbb{Q}_{\theta, \gamma, \zeta }$, a real number can be mapped to the closest point in $\mathbb{Q}_{\theta, \gamma, \zeta}$. To limit the quantization error, the mapping $m(\cdot)$ needs to satisfy 
\begin{equation}
    |m(\varphi)-\varphi| \leq \theta^{-\zeta }, \forall \varphi \in\left[-\theta^{\gamma}, \theta^{\gamma}\right]\label{28ss}
\end{equation}
where the quantization error is restricted by the resolution within the range of $\mathbb{Q}_{\theta, \gamma, \zeta }$. To map the real-number set onto the integer set $\mathbb{Z}$, we simply scale $\mathbb{Q}_{\theta, \gamma, \zeta }$ by $\theta^\zeta $ as 
\begin{equation}
  \mathbb{Z}_{\theta, \gamma, \zeta } = \theta^{\zeta } \mathbb{Q}_{\theta, \gamma, \zeta } {=} \left\{-\theta^{\gamma+\zeta },-\theta^{\gamma+\zeta }{+}1, \ldots, \theta^{\gamma+\zeta }{-}1\right\} \label{29ss}
\end{equation}
where $\mathbb{Z}_{\theta, \gamma, \zeta } \subseteq \mathbb{Z}$ denotes the fixed-point set in the integer field. Moreover, the SS requires the inputs to be within the field $\mathbb{E}$. Therefore, we further map each element in $z \in \mathbb{Z}_{\theta, \gamma, \zeta }$ onto $\mathbb{E}$ with the modular operation as 
\begin{equation}
    g(z) = z \bmod e. 
    \label{30ss}
\end{equation}
Note that $z \in \mathbb{Z}_{\theta, \gamma, \zeta }$ can be any negative integer, and the modular operation in \eqref{30ss} will change the sign of a negative input, i.e., $g(\hat{z}) = \hat{z}+e$ for $\hat{z}<0$. To address the negative integer operation, we introduce the partial inverse of $g(\cdot)$ as
\begin{equation}
  \psi(z)=\left\{\begin{array}{ll}z-e & \text { if } z \geq \frac{e}{2}, \\ z & \text { otherwise. }\end{array}\right. \label{31ss}
\end{equation}
Therefore, we can readily obtain $z=\psi(g(z)),\forall z\in \mathbb{E}$.

\subsection{SS-based Privacy-Preserving Algorithm}
\subsubsection{Shamir's secret sharing scheme} 
Before introducing the privacy-preserving algorithm design, we first briefly introduce Shamir's SS scheme \cite{abb_shamir1979share} which merits an efficient and  lightweight private information distribution structure. Suppose a manager (secret holder) seeks to distribute a secret $\omega$ to specific agents and mandates the cooperation of at least $d$ agents to retrieve the secret. In such needs, Shamir's SS is grounded on the following idea of Lagrange interpolation for secret distribution and recovery.

\noindent \textbf{Theorem 1} (\textit{Polynomial interpolation}\cite{abb_humpherys2020foundations})\textbf{.} Let $\{(\varsigma_{1}, y_{1}), \ldots,$ $(\varsigma_{d}, y_{d})\}$ $\subseteq \mathbb{R}^{2}$ be a set of points whose  values of $\varsigma_\imath$ are all distinct. Then there exists a unique polynomial $\mathcal{Y}$ of degree $d-1$ that satisfies $y_\imath = \mathcal{Y}(\varsigma_\imath), \forall \imath=1,\ldots,d$.  \hfill $\blacksquare$

In SS-based schemes, the manager first constructs a random polynomial of degree $d-1$ as 
\begin{equation}
        y(z) = \omega + a_1z + \cdots + a_{d-1}z^{d-1} \label{f(z)}
    \end{equation}
where $\omega$ denotes an integer secret, $a_1,\ldots, a_{d-1}$ are random coefficients that are uniformly distributed in the field $\mathbb{E} \triangleq [0,e)$, and $e$ denotes a prime number that is larger than $\omega$. Secondly, the manager calculates the outputs of \eqref{f(z)} with non-zero integer inputs, e.g., setting $\tau=1,\ldots,n$ to retrieve $(\tau,y(\tau))$ where $ y_\tau^\Pi = y(\tau) \bmod e$. Then, the share $y_\tau^\Pi$ is distributed to agent $\tau$. Lastly, at least $d$ agents with shares  are required to reconstruct the polynomial based on \textbf{Theorem 1} and hence recover the secret $\omega$ by 
\vspace{-1mm}
\begin{equation}
   \omega =  \sum_{\tau=1}^{d}y_\tau^\Pi \prod_{\upsilon =0  \atop \upsilon \neq \tau}^{d} \frac{\upsilon}{\upsilon-\tau}.
   \label{recosntruction}
\end{equation}

\vspace{-1mm}
\subsubsection{Proposed  privacy-preserving algorithm}

\begin{figure}[!htb]
\vspace{-1mm}
    \centering
\includegraphics[width=0.45\textwidth, trim={0.2cm 0cm 0.2cm 0.22cm},clip]{ 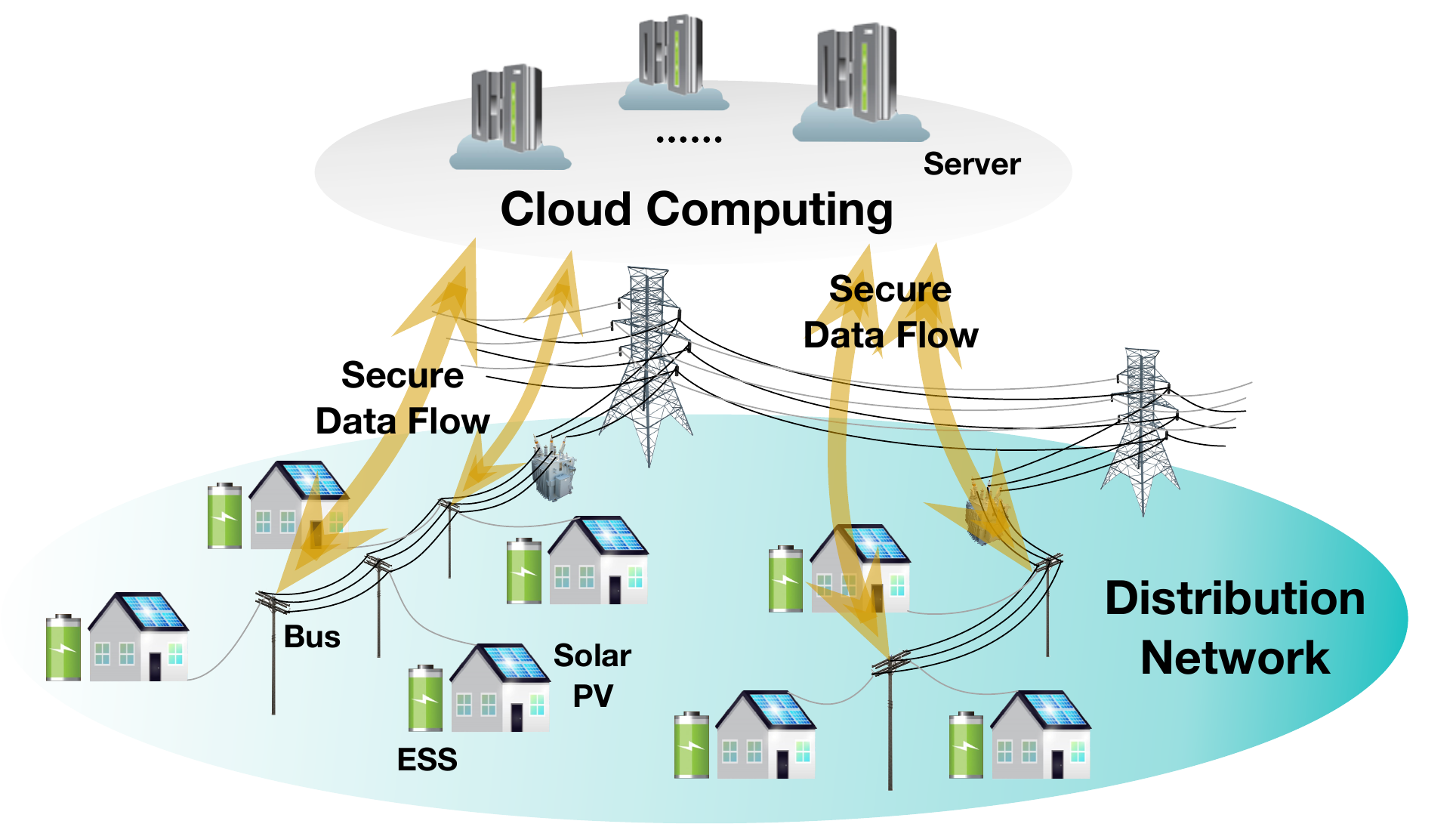}
    \caption{Two-layer privacy-preserving computing structure for DER control in distribution networks.}
    \label{fig2}
     \vspace{-2mm}
\end{figure}

We next present the proposed two-layer decentralized privacy-preserving algorithm based on SS in a bus-level aggregation and control architecture, to achieve privacy preservation and scalability concurrently. In the distribution network layer, all DERs' decision variables are updated in parallel, and only masked data are sent from each bus to the servers. In the cloud computing layer, the servers calculate the aggregated messages and  distribute them to the related buses. The computing structure of the proposed privacy-preserving algorithm is shown in Fig. \ref{fig2}.


Let $\mathbb{C}$ denote the set of clouds and $c \geq 2$ denotes the total number of clouds. The $i$th bus generates a random polynomial of order $d-1$ using \eqref{f(z)} to obtain  
\begin{equation}
y_i^{(\ell)}(z) = \omega_i^{(\ell)} + a_{i,1}^{(\ell)}z + \cdots + a_{i,d-1}^{(\ell)}z^{d-1} \label{34s}
\end{equation}
where $2 \leq d \leq c$, $\omega_i^{(\ell)}$ denotes the secret of bus $i$ at the $\ell$th iteration, $\ell$ denotes the iteration number, and $a_{i,1}^{(\ell)},\ldots, a_{i,d-1}^{(\ell)}$ denote random coefficients that are uniformly distributed in the field $\mathbb{E}$. Note that for a vector secret such as $\bm{p}_i$, we refer to an elementwise calculation of the vector using \eqref{34s} by default.

At the $\ell$th iteration, the $u$th cloud firstly generates a random integer $\alpha_u^{(\ell)}$, then it broadcasts $\alpha_u^{(\ell)}$ to all the buses. Subsequently, the $i$th bus can calculate $y_i^{(\ell)}(\alpha_u^{(\ell)})$, $\forall u=1,\ldots,c$ using the received inputs based on \eqref{34s}. Finally, the $i$th bus sends $y_i^{(\ell)}(\alpha_u^{(\ell)})$ back to the $u$th cloud. 
Note that the coupling term $\bm{\pi}_i \Tilde{\bm{P}}$ in \eqref{27sss} is a linear combination of all $\bm{p}_i$'s that requires the private generation/consumption details from the buses. Therefore, 
a secure computation framework of $\bm{\pi}_i \Tilde{\bm{P}}$ is required to preserve the  privacy of buses and DER owners.


Suppose the clouds are aware of the network topology matrix $\bm{Z}$ which contains no private information of the buses or DERs. In order to calculate the aggregated information $\bm{\pi}_i \Tilde{\bm{P}}$ for bus $i$, the $u$th cloud firstly  multiplies the received outputs $y_{1}(\alpha_u^{(\ell)}),\ldots ,y_{n}(\alpha_u^{(\ell)})$  utilizing the coefficients of  $\bm{\pi}_i$ to obtain
\begin{equation}
   \{
    \alpha_u^{(\ell)},  \bm{\pi}_i(1)y_{1}^{(\ell)}(\alpha_u^{(\ell)}),
    \ldots, \bm{\pi}_i(n)y_{n}^{(\ell)}(\alpha_u^{(\ell)})\}
    \label{36ss}
\end{equation}
Then, the $u$th cloud  sums the outputs in \eqref{36ss} to obtain a new pair of input and output as 
\begin{equation}
  \Bar{\mathcal{A}}_{u,i} =\{\alpha_u^{(\ell)},
\sum_{\hat{\imath}=1}^{n}\bm{\pi}_i(\hat{\imath}) ~ y_{\hat{\imath}}^{(\ell)}(\alpha_u^{(\ell)})\}. \label{38ssss}
\end{equation}
Finally, the $u$th cloud calculates $\Bar{\mathcal{A}}_{u,i}$, $\forall i=1,\ldots,n$ and broadcasts the new input-output share $\Bar{\mathcal{A}}_{u,i}$ to the $i$th bus.

Therefore, after receiving new shares from in total $c$ clouds servers, the $i$th bus now has access to 
\begin{equation}
    \tilde{\mathcal{A}}_{i} = \left\{
    \alpha_{\hat{\jmath}}^{(\ell)}, \sum_{\hat{\imath}=1}^{n}\bm{\pi}_i(\hat{\imath}) ~ y_{\hat{\imath}}^{(\ell)}(\alpha_{\hat{\jmath}}^{(\ell)}), \forall \hat{\jmath} = 1, \ldots, c\right\}.
    \label{37ss}
\end{equation} Note that $\tilde{\mathcal{A}}_{i}$ contains in total $c$ shares that can construct a new polynomial of the form
\begin{equation}
\tilde{y}_i^{(\ell)}(z) = \bm{\pi}_i \tilde{\bm{P}} + \tilde{a}_{i,1}^{(\ell)}z + \cdots + \tilde{a}_{i,d-1}^{(\ell)}z^{d-1} \label{43}
\end{equation}
whose constant term is exactly $\bm{\pi}_i \Tilde{\bm{P}}$.

\begin{figure}[!htb]
    \centering
    \includegraphics[width=0.4\textwidth, trim={0cm 0.1cm 0.2cm 0.1cm},clip]{ 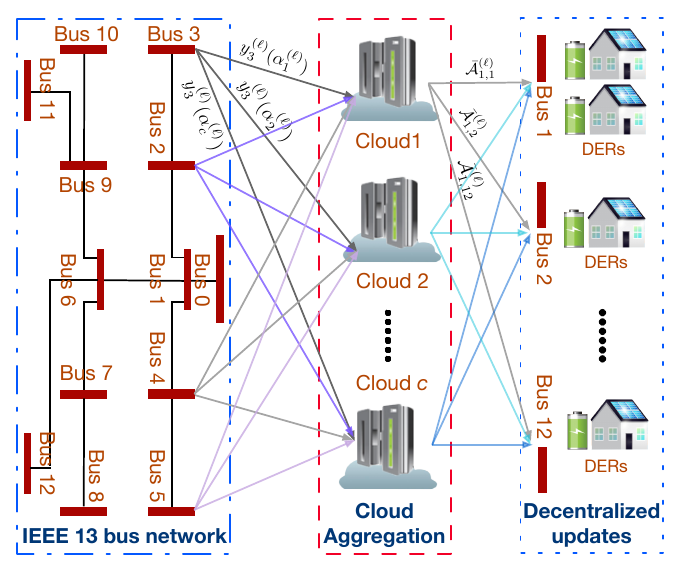}
  \caption{Information exchange structure between the distribution network and cloud servers (only the messages sent from bus 3 and cloud 1 are labeled).}
    \label{information_exchange}
\end{figure}

During this information exchange process, each bus only sends a single share to each server so that a single cloud server is incapable of reconstructing the secret based on the received shares, and herein cannot infer agents' true decision variables. The cloud servers only need to calculate aggregated messages using   outputs of randomized polynomials. The details of the proposed method are presented via \textbf{Algorithm \ref{alg_1}}.

\begin{algorithm}
\caption{Decentralized SS-based  privacy-preserving DER control strategy}
\begin{algorithmic}[1]

\State Agents initialize decision variables, tolerance $\epsilon_0$,  basis $\theta$, magnitude $\gamma$, resolution $\zeta $, iteration counter $\ell=0$, and maximum iteration $\ell_{max}$.

\While{ $\epsilon_{\nu(\sigma)}^{(\ell)} > \epsilon_0$ and $\ell < \ell_{max}$}

\State Each bus performs real number to integer transformation using \eqref{27ss}-\eqref{30ss}, then obtains the integer secret $\omega_i^
{(\ell)}$.

\State The $u$th cloud generates a random integer $\alpha_u^
{(\ell)}$, then broadcasts $\alpha_u^
{(\ell)}$ to all the buses.

\State The $i$th bus generates a random polynomial $y_i^
{(\ell)}(z)$ using \eqref{34s} with $\omega_i^
{(\ell)}$ as the constant term, calculates the outputs using  $\alpha_1^
{(\ell)},$ $\ldots,$ $\alpha_c^
{(\ell)}$ to obtain $y_i^
{(\ell)}(\alpha_1^
{(\ell)}),$ $\ldots,$ $y_i^
{(\ell)}(\alpha_c^
{(\ell)})$, then sends $y_i^
{(\ell)}(\alpha_u^
{(\ell)})$ to the $u$th cloud.

\State  The $u$th cloud formulates $\Bar{\mathcal{A}}_{u,i}$ in \eqref{38ssss}, then broadcasts $\Bar{\mathcal{A}}_{u,i}$ to the $i$th bus.


\State The $i$th bus formulates $\tilde{\mathcal{A}}_{i}$ in \eqref{37ss},  reconstructs the aggregated secrets using $c$ shares to obtain $\bm{\pi}_{i}\Tilde{\bm{P}}$, then calculates $\mathcal{C}_p$ in \eqref{27sss}.

\State  The $i$th bus transforms $\mathcal{C}_p$ back to real numbers using \eqref{31ss}, then decision variables $\tilde{\bm{p}}_{\nu}^{(\ell)}$ or $\hat{\bm{p}}_{\sigma}^{(\ell)}$ of DERs connected at bus $i$ are updated by PGM using \eqref{eq1}. The $i$th bus calculates the error $\epsilon_{\nu}^{(\ell)}$ or $\epsilon_{\sigma}^{(\ell)}$.

\State  $\ell=\ell+1$. 

\EndWhile
\end{algorithmic}
\label{alg_1}
\end{algorithm}

\textbf{Algorithm \ref{alg_1}} can achieve privacy preservation while maintaining exact solutions as non-privacy PGM-based methods. The decision variables will be continuously updated till the convergence errors $\epsilon_{\nu}^{(\ell)}\triangleq \|\tilde{\bm{p}}_{\nu}^{(\ell)}-\tilde{\bm{p}}_{\nu}^{(\ell-1)}\|_2^2$ and $\epsilon_{\sigma}^{(\ell)} \triangleq  \| \hat{\bm{p}}_{\sigma}^{(\ell)} - \hat{\bm{p}}_{\sigma}^{(\ell-1)} \|_2^2$  are smaller than the threshold $\epsilon_0$. The correctness of \textbf{Algorithm \ref{alg_1}} is presented via \textbf{Theorem 2}.

\noindent \textbf{Theorem 2}  (\emph{Correctness})\textbf{.} Let $\mathbb{E}$ denote the domain of the input secrets $\omega_{1}, \ldots, \omega_{n}$, and $\mathcal{C}_{p}$  denote the desired outputs. Then, \textbf{Algorithm \ref{alg_1}} satisfies:
\begin{equation}
\operatorname{Pr}\left[\forall c \geq d, \operatorname{Rec}\left(\mathbb{A}, \mathbb{E}, \bm{Z}, \bar{\delta}_1,  \theta, \gamma, \zeta  \right)= \mathcal{C}_p \right]=1
\end{equation}
where $\mathbb{A} = \{\tilde{\mathcal{A}}_1, \ldots, \tilde{\mathcal{A}}_c\}$ denotes the set of shares from agents, $\operatorname{Pr}[\cdot]$ denotes probability, and $\operatorname{Rec}(\cdot)$ denotes the secret reconstruction operation. \hfill $\blacksquare$

\textbf{Theorem 2} states that \textbf{Algorithm \ref{alg_1}
}can correctly retrieve the aggregated information $\mathcal{C}_p$ which would be further used to achieve exact primal and dual updates. 

The detailed proof of \textbf{Theorem 2} can be found in \textsc{Appendix} \ref{Theorem2_proof}.

\noindent \textbf{Remark 1 :} Though \textbf{Algorithm \ref{alg_1}} is  developed based on bus-level aggregation and control, it can also be extended to the DER-level aggregation and control. In DER-level aggregation and control, each DER is required to generate a polynomial in \eqref{34s} and act as an independent agent in secret reconstruction using \eqref{37ss}. Besides, depending on the practical applications, DERs can also be clustered and controlled  by the household or district where the new clusters act as agents, following the similar design of  \textbf{Algorithm \ref{alg_1}}.  \hfill $\square$

\noindent \textbf{Remark 2:}
The multi-server architecture seamlessly integrates the SS scheme into DER aggregation and control. Shares generated from buses were  aggregated and broadcasted to the buses by a group of servers for the purpose of secret retrieval. The aggregation task is distributed to multiple servers to ensure that a single server cannot retrieve any secrets. \hfill $\square$

\subsection{Privacy Analysis}

The proposed approach aims at protecting the decision variables of the DERs whose disclosure can lead to the leakage of customers' sensitive information. To resolve this issue, \textbf{Algorithm \ref{alg_1}} achieves privacy preservation against two types of adversaries, including \emph{honest-but-curious-agent} who follows the algorithm but may utilize the possessed and received data to infer the private information of other agents, and \emph{external eavesdroppers} who wiretap and intercept exchanged messages from communication channels.

\noindent \textbf{Proposition 1:} (\textit{Secure cloud computing})\textbf{.} In \textbf{Algorithm \ref{alg_1}}, any cloud number less than $d-1$ cannot infer any information of the aggregated decision variables $\mathcal{C}_p$. \hfill $\blacksquare$

\textbf{Proposition 1} presents the security of the proposed algorithm against corrupted clouds. Based on the polynomial interpolation in  
\textbf{Theorem 1}, at least $d$ clouds are required to retrieve any secret through collusion. 


\textbf{Proposition 1} is proved based on the correctness analysis. Please refer to \textsc{Appendix} \ref{Proposition_1_proof} for the detailed proof.

\noindent \textbf{Assumption 1.} At least one communication link of an individual agent is secure against external eavesdroppers. \hfill $\blacksquare$

\textbf{Assumption 1} is essential and generically used in SS-based schemes. Given $d$ pairs of shares sent via different communication links, i.e.,  $\{(\varsigma_{1}, y_{1}), \ldots,$ $(\varsigma_{d}, y_{d})\}$ $\subseteq \mathbb{R}^{2}$, if an external eavesdropper wiretap all communication links to gain access to the shares, then it can simply deduce the secret by Lagrangian interpolation using \textbf{Theorem 1}.

\noindent \textbf{Theorem 3} (\textit{Privacy preservation against adversaries})\textbf{.} By using \textbf{Algorithm \ref{alg_1}}, the following two statements stand:
\begin{enumerate}
\item \textbf{Algorithm \ref{alg_1}} securely computes and updates the decision variables between agents in the presence of honest-but-curious agents.   
    \item External eavesdroppers learn no private information of the agents. \hfill $\blacksquare$
\end{enumerate}

\textbf{Theorem 3} gives privacy preservation guarantees in the presence of honest-but-curious agents and external eavesdroppers.
The privacy preservation of \textbf{Algorithm \ref{alg_1}} can be proved from secure multi-party computation (SMC) perspective. Before giving detailed privacy analyses and proofs, we first introduce some concepts of SMC. 

\noindent \textbf{Definition 1} (\emph{Computational indistinguishability}\cite{abb_goldreich2009foundations})\textbf{.}
Let $\{D_{\varkappa}\}_{{\varkappa\in {\mathbb{N}}}}$ and $\{E_{\varkappa}\}_{{\varkappa\in {\mathbb{N}}}}$ be two distribution ensembles with security parameter $\varkappa$; If for any non-uniform probabilistic polynomial-time algorithm $\mathcal{G}$, $\delta(\varkappa)$ is negligible, where
\begin{equation}
\delta(\varkappa)=\left|\underset{x_1 \leftarrow D_{\varkappa}}{\operatorname{Pr}}[\mathcal{G}(x_1)=1]-\underset{x_2 \leftarrow E_{\varkappa}}{\operatorname{Pr}}[\mathcal{G}(x_2)=1]\right|
\end{equation}
we say that $\{D_{\varkappa}\}_{{\varkappa\in {\mathbb  {N}}}}$ and $\{E_{\varkappa}\}_{{\varkappa\in {\mathbb  {N}}}}$ are computationally indistinguishable, denoted as $D_{\varkappa} \stackrel{c}{\equiv} E_{\varkappa}$.  \hfill $\blacksquare$

Therefore, \textbf{Definition 1} states that any polynomial-time algorithm cannot distinguish two computationally indistinguishable ensembles because the outputs of those algorithms do not significantly differ. In what follows, \textbf{Definition 2} presents the standard privacy notion in SMC.

\noindent \textbf{Definition 2} (\hspace{1sp}\cite{abb_evans2018pragmatic,goldreich1998secure})\textbf{.} Let $\Pi$ be an $m$-party protocol for computing the outputs of function $\mathcal{F}(\bar{x})$ where $\bar{x}=\{x_{1}, \ldots, x_{m}\}$ and $\mathcal{F}_{\rho} (\bar{x})$ denotes the $\rho$th output of $\mathcal{F} (\bar{x})$. Let $\mathbb{M}=\left\{M_{1}, \ldots, M_{m}\right\}$ denote the set of parties. The view of the $\rho$th party during the execution of  $\Pi$ is denoted by $\operatorname{VIEW}_{\rho}^{\Pi}(\bar{x})$. We say that $\Pi$ privately computes $\mathcal{F}(\bar{x})$ if there exists a polynomial-time algorithm $\mathcal{S}$, such that for every party $M_\rho$ in $\mathbb{M}$, we have
\begin{equation}
\mathcal{S}(\rho, x_{\rho}, \mathcal{F}_{\rho}(\bar{x})) \stackrel{c}{\equiv} \operatorname{VIEW}_{\rho}^{\Pi}(\bar{x}). 
\label{42ss}
\end{equation} \hfill $\blacksquare$

\textbf{Definition 2} states that the security of an $m$-party protocol can be evaluated based on computational indistinguishability, i.e., the view of the parties can be efficiently simulated based solely on their inputs and outputs. In other words, SMC allows a group of participants to learn the correct outputs of some agreed-upon function applied to their private inputs without revealing anything else. The 
theoretical underpinnings of \textbf{Definition 1} and  \textbf{Definition 2} can help prove that \textbf{Algorithm \ref{alg_1}} securely computes $\bm{\pi}_1 \tilde{\bm{P}},\ldots,\bm{\pi}_n \tilde{\bm{P}}$ between the agents. 

The detailed proofs of \textbf{Theorem 3} can be found in \textsc{Appendix} \ref{Theorem_3_proof}.

\begin{figure*}[!t]
 \vspace{-4mm}
  \centering
  \subfigure[Heterogeneous baseline loads of 24 houses]{%
    \includegraphics[width=0.24\textwidth]{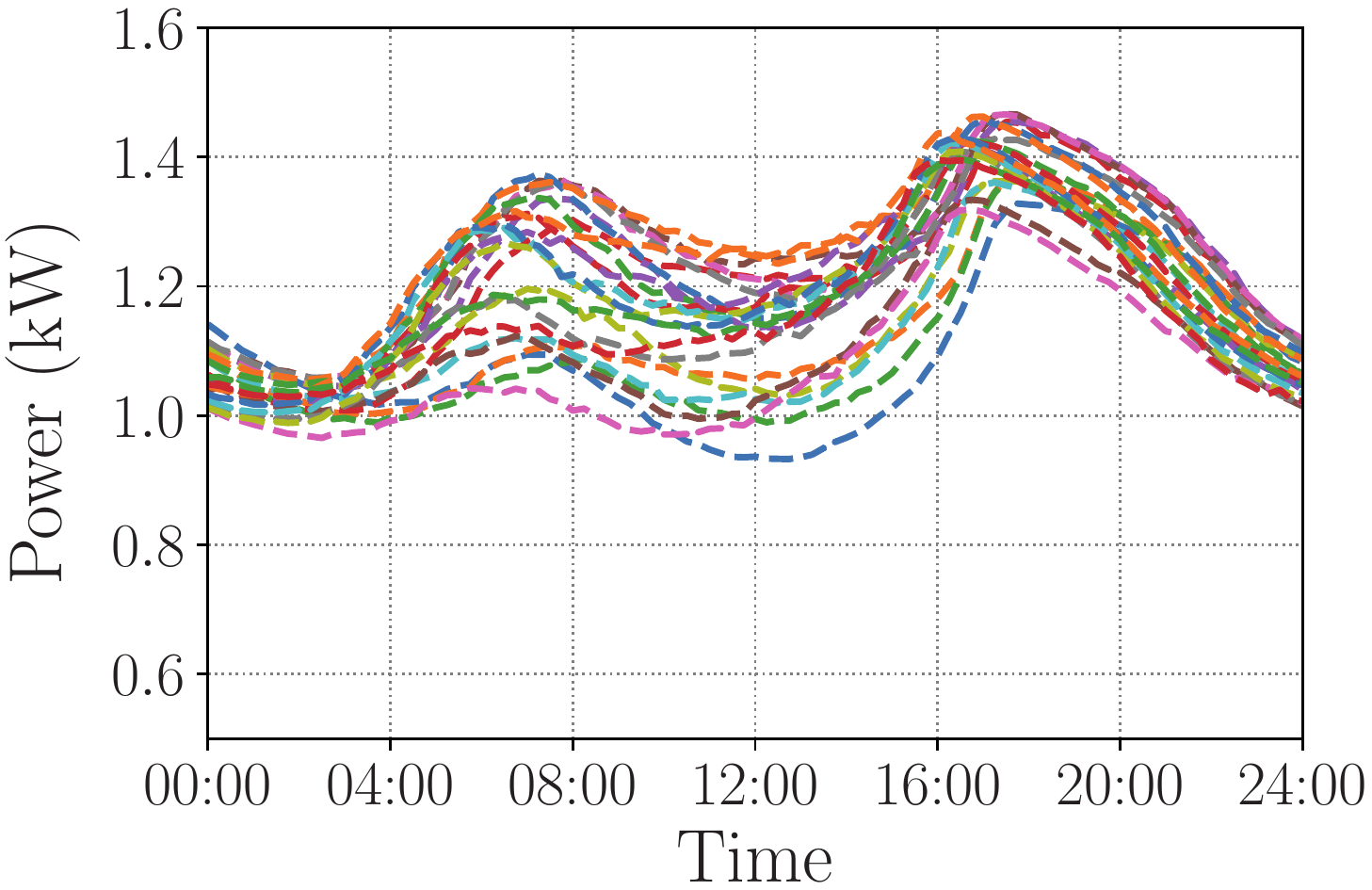}
    \label{fig4_a}%
  }%
  \hfill
  \subfigure[Solar power injection of  24 houses]{%
    \includegraphics[width=0.24\textwidth]{ 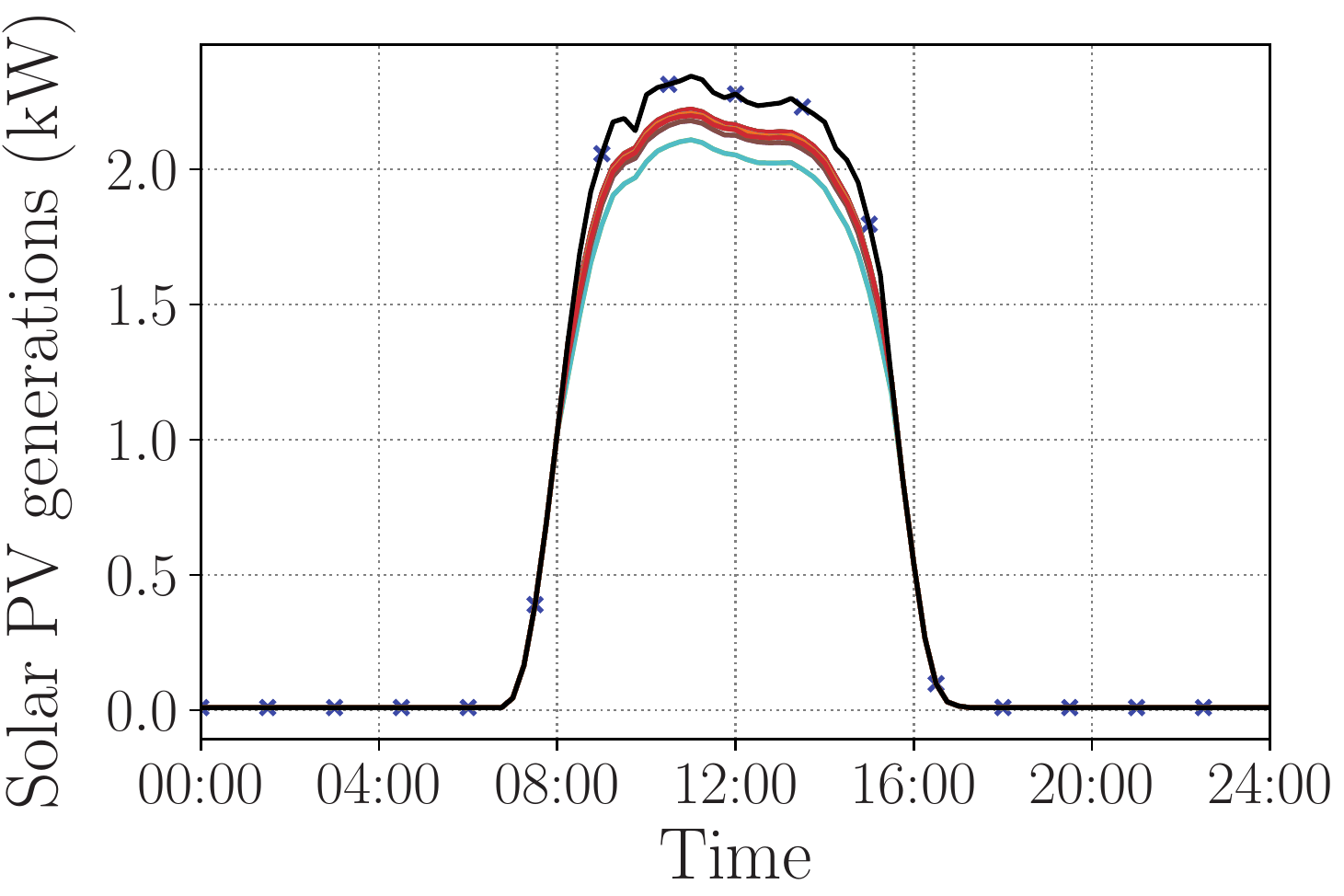}%
    \label{fig4_b}%
  }%
  \hfill
  \subfigure[Charging and discharging power from 24 ESSs]{%
    \includegraphics[width=0.255\textwidth, trim={0cm 0cm 0cm 0cm},clip]{ 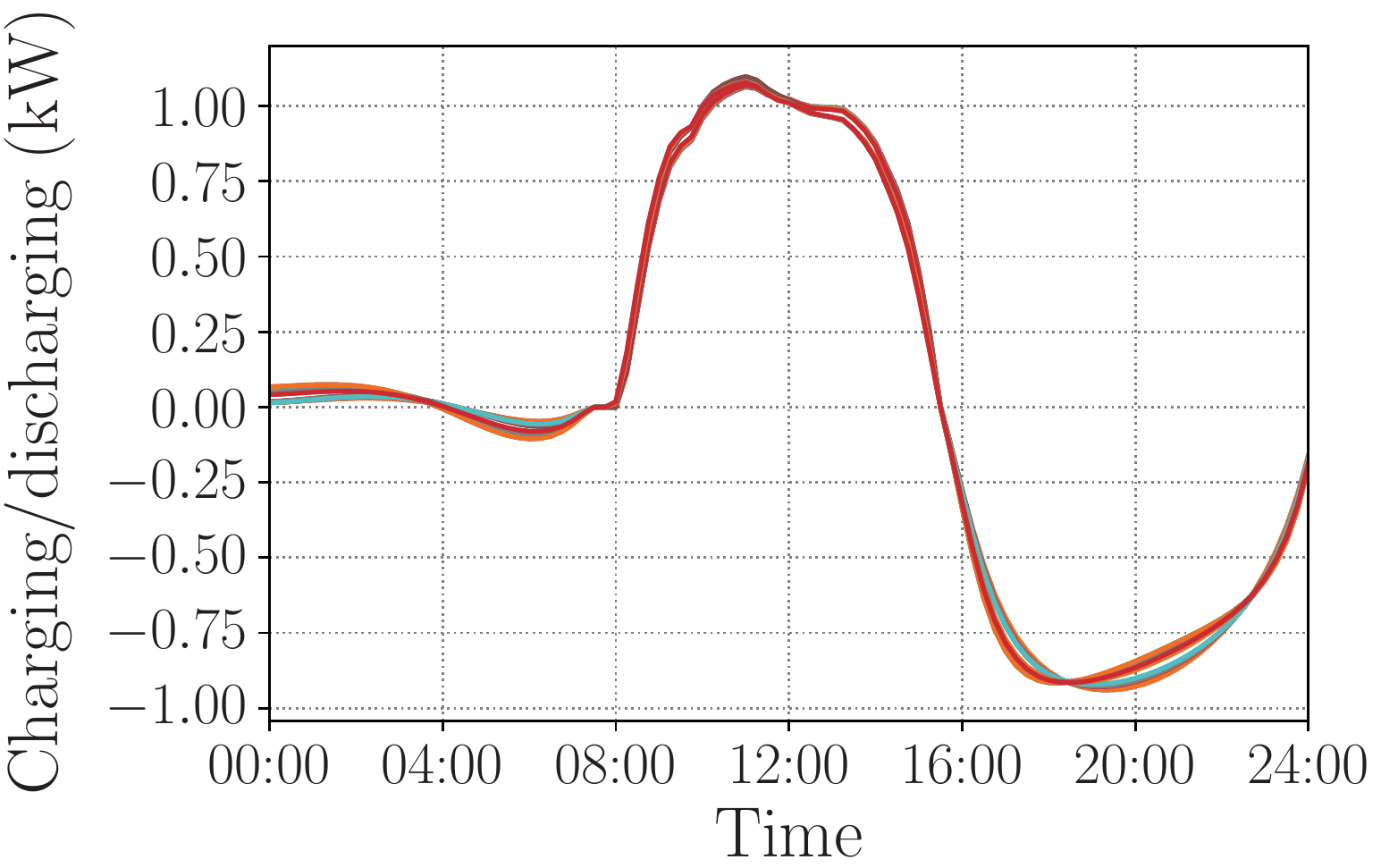}
    \label{fig4_c}%
  }%
  \hfill
  \subfigure[Power flows of 12 lines in the distribution network]{%
    \includegraphics[width=0.236\textwidth, trim={0cm 0.05cm 0cm 0cm},clip]{ 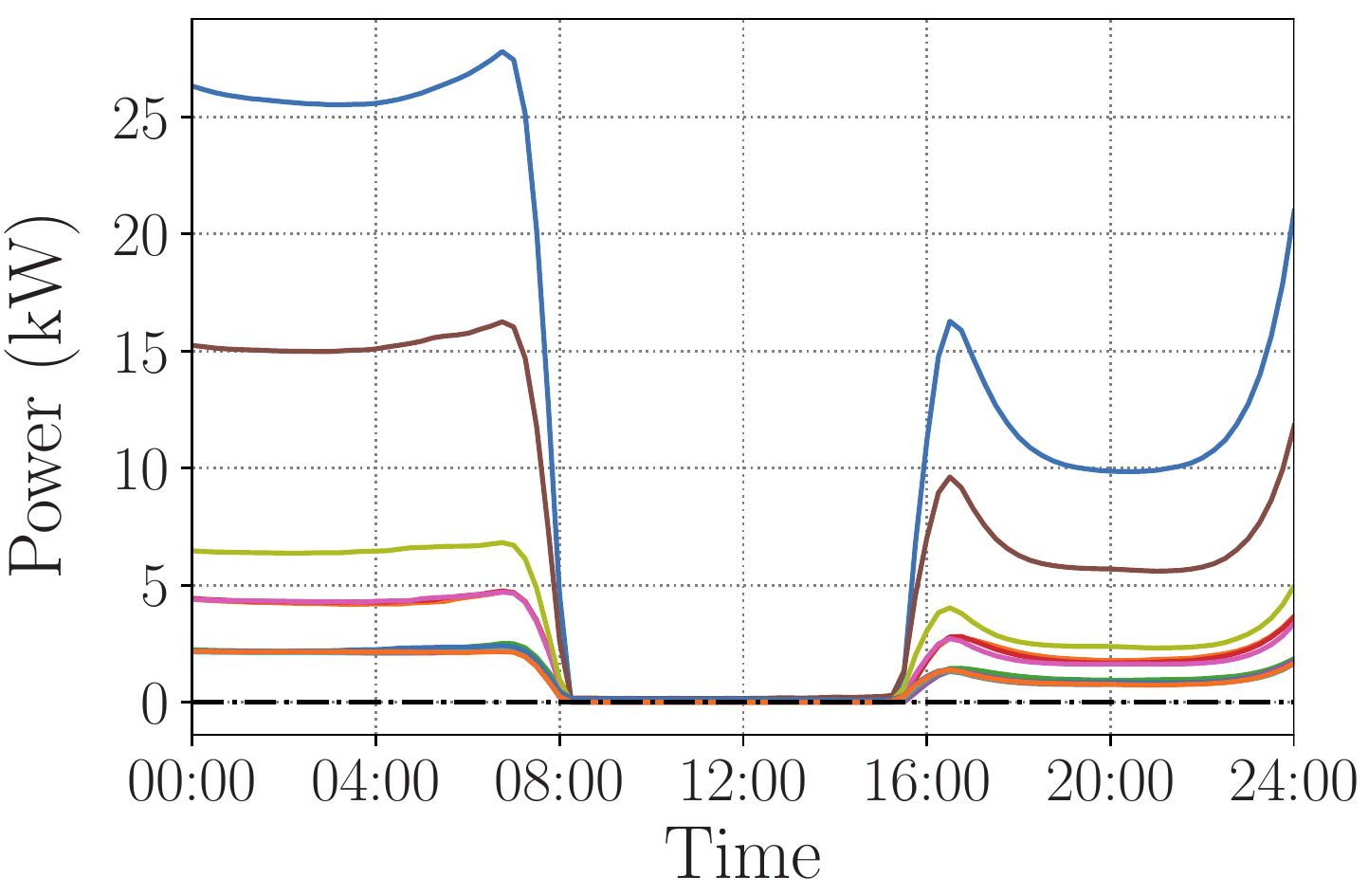}%
    \label{fig4_d}%
  }%
  \vspace{-2mm}
  \caption{The optimal solutions of \eqref{p2} by controlling DERs in the distribution network.}
  \label{fig4}
    \vspace{-4mm}
\end{figure*}


\section{Simulation Results}
\label{Simulation_results}

A simplified single-phase IEEE 13-bus test feeder \cite{abb_liu2017decentralized} is used to verify the proposed decentralized privacy-preserving DER control strategy. In specific, each bus, except the feeder head, is assumed to be connected with 2 houses and each house is equipped with an ESS and 5 solar panels that can generate maximum 2.5 kW solar output. The maximum capacity of all residential ESSs are 10 kWh, the initial SoCs of all ESSs are uniformly set to be $4$ kWh, and the maximum charging and discharging rates are $\pm3$ kW, respectively \cite{ESS}. The forecasted solar PV generation is chosen from 01/01/2021 with $\Delta T = 15$ mins in California from CAISO \cite{CISO}. 

In total $c=4$ clouds are responsible for message  aggregation and distribution. The degree of all polynomials is set to be $d-1=3$ and the integer field is chosen as $\mathbb{E} = [0, 2^{31}-1)$. For the fixed-point number quantization, the basis, magnitude, and resolution are uniformly set to be $\theta=2$, $\gamma=27$, and $\zeta=4$, respectively. For the distribution network shown in Fig. \ref{13bus}, all 24 houses are assumed to be located in the same area with identical solar radiation. The baseline load profiles of all houses are shown in Fig. \ref{fig4_a} \cite{CISO}. The primal and dual step sizes are chosen based on experience to be $\alpha_{\nu,\ell}^{v} =2.3$, $\alpha_{\sigma,\ell}^{e} = 1.8$, and $\beta_{\mu_{l\iota},\ell} = 5\times 10^{-4}$, respectively. Note that only the lower bound of power flow limits in \eqref{17sss} is active, herein, only the results related to $\bm{\mu}_{l\iota}$ are presented. 

\begin{figure}[!htb]
 \vspace{-2mm}
  \centering
  \subfigure[\scriptsize Convergence of solar PVs' decision variables $\tilde{\mathbf{p}}_{\nu}$]{%
    \includegraphics[width=0.235\textwidth,trim={0.1cm 0cm 0cm 0.25cm},clip]{ 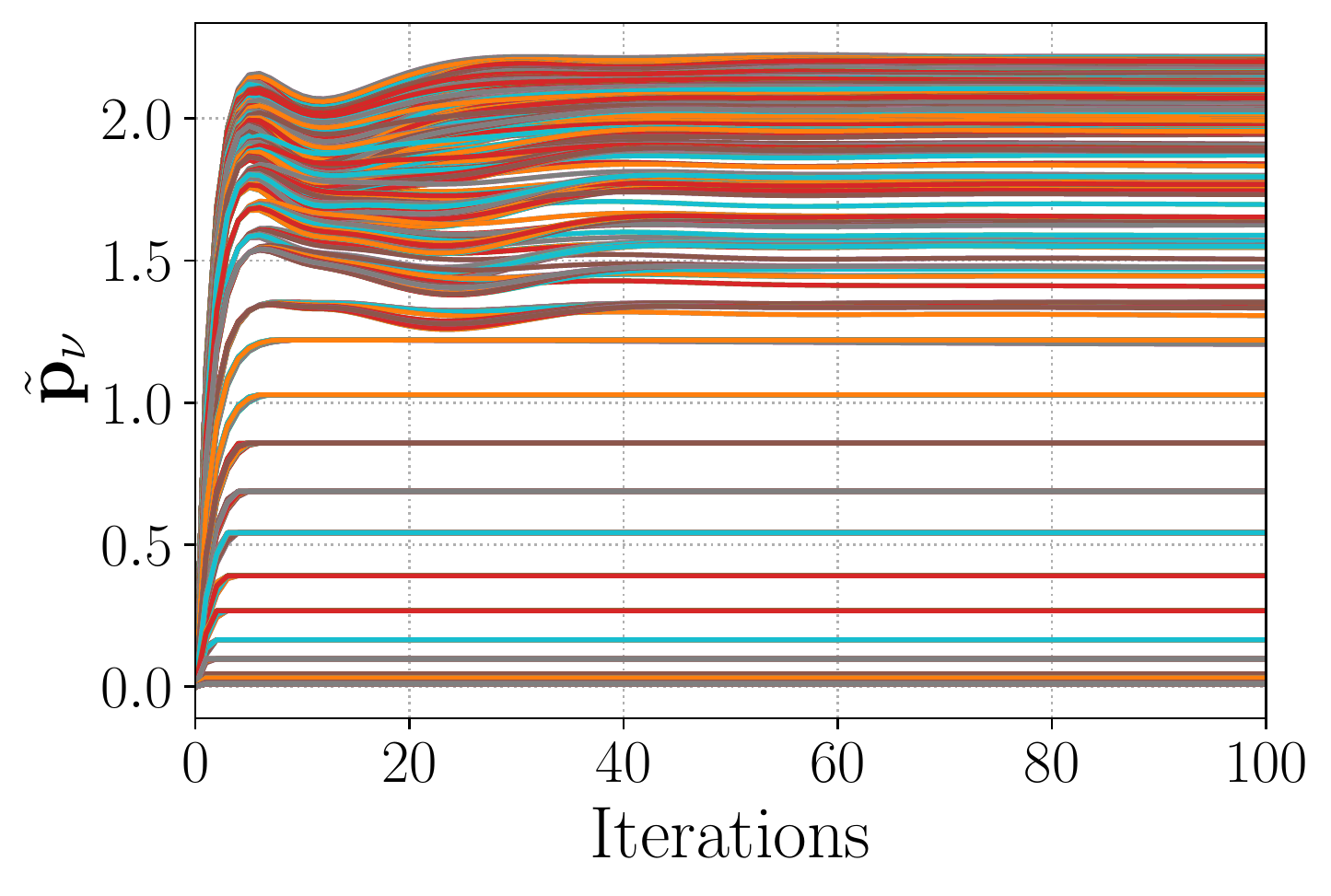}%
    \label{fig_6a}%
  }%
  \hfill
  \subfigure[\scriptsize
  Convergence of the dual variable $\bm{\mu}_{l\iota}$ ]{%
    \includegraphics[width=0.24\textwidth,trim={0.1cm 0cm 0cm 0.25cm},clip]{ 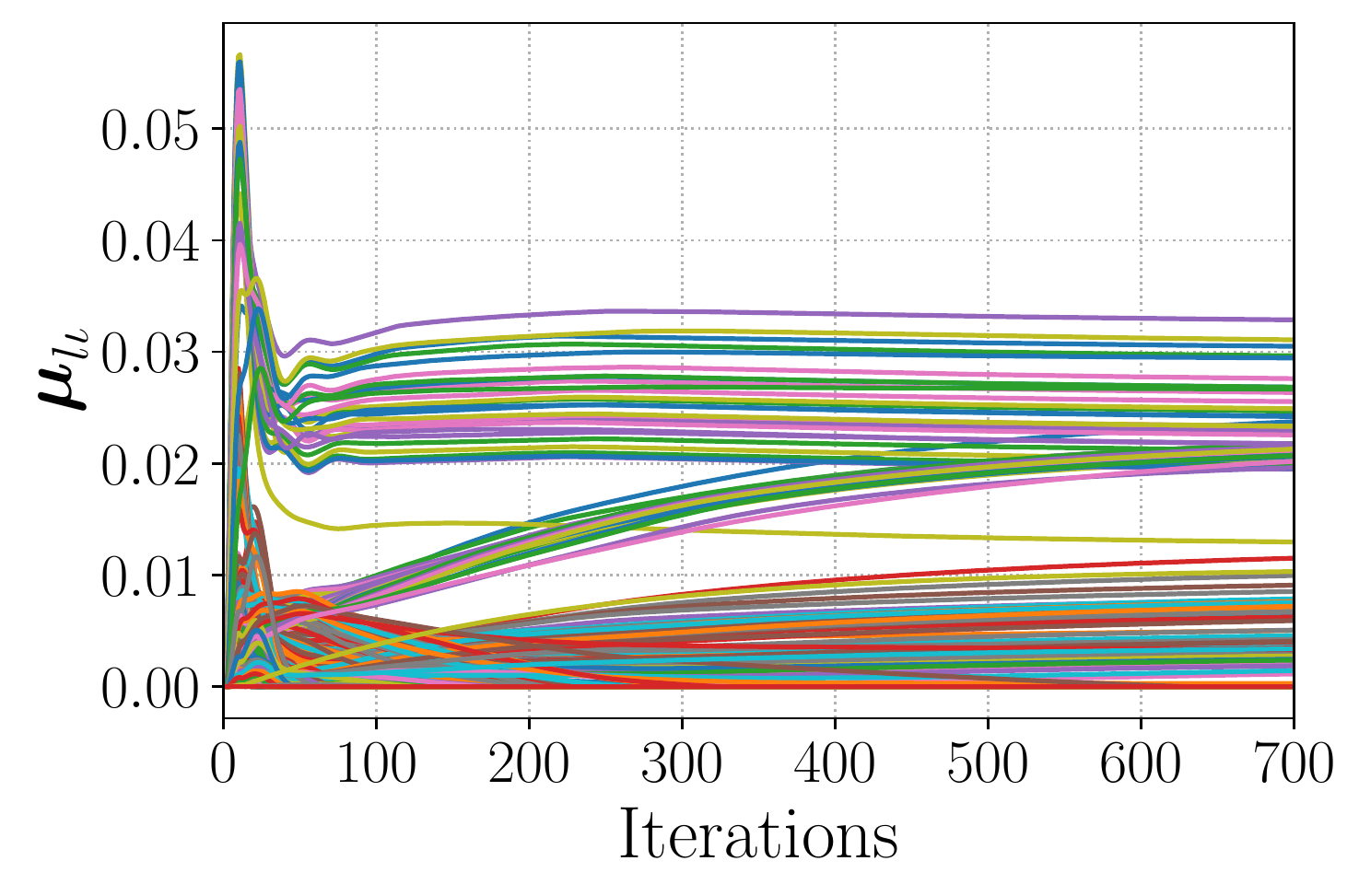}%
    \label{fig_6b}%
  }%
\vspace{-2mm}
  \caption{Convergence of the primal and dual variables}
  \label{fig6}
\end{figure}

Fig. \ref{fig4_b} and Fig. \ref{fig4_c} show the  active power generations and the charging/discharging  power from the solar PVs and ESSs, respectively. At around 12:00, the solar PVs generate the maximum amount of energy, and the ESSs charge at peak rates. After 16:00, energy stored in ESSs is extracted to supply in-home use and compensate for the power loss in the distribution network. The power flows of 12 lines are shown in Fig. \ref{fig4_d} where no inverse flows occur. Moreover, accurate primal and dual solutions are achieved without affecting the anticipated primal-dual convergence. The iterative solutions of the primal and dual variables are shown in Fig. \ref{fig6}. 
\begin{figure}[!htb]
\vspace{-5mm}
    \centering
\includegraphics[width=0.4\textwidth, trim={0cm 0cm 0cm 0cm},clip]{ 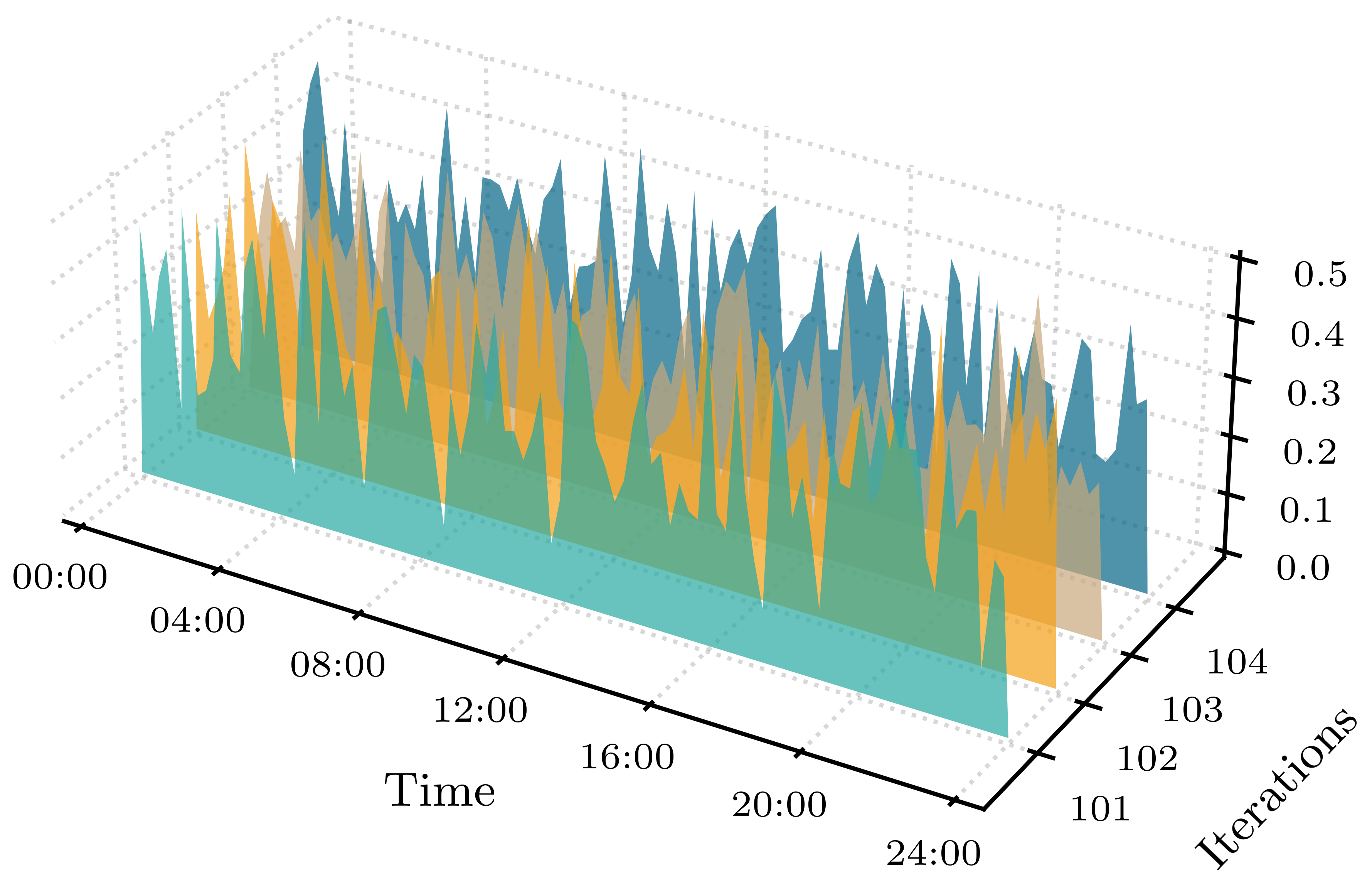} 
\vspace{-2mm}
\caption{Random shares generated by Bus 6 at different iterations}
    \label{fig7}
    \vspace{-3mm}\end{figure}
Fig. \ref{fig7} presents normalized shares generated by Bus 6 using the random polynomial $y_6^
{(\ell)}(z)= \omega_6^{(\ell)} + a_1^
{(\ell)}z + a_{2}z^{2} + a_{3}^
{(\ell)}z^{3}$ where the coefficients $a_i^
{(\ell)}$, $i=1,2,3$ are randomized at each iteration and different time slots. The privacy preservation of \textbf{Algorithm \ref{alg_1}} against external eavesdroppers are guaranteed because external eavesdroppers have insufficient information in polynomial reconstruction by wiretapping the transmitted shares. Without loss of generality, suppose bus 6 is honest-but-curious. Fig. \ref{fig8} \begin{figure}[!htb]
     \vspace{-3mm}
    \centering
\includegraphics[width=0.5\textwidth, trim={0.5cm 0.7cm 0.5cm 0.5cm},clip]{ 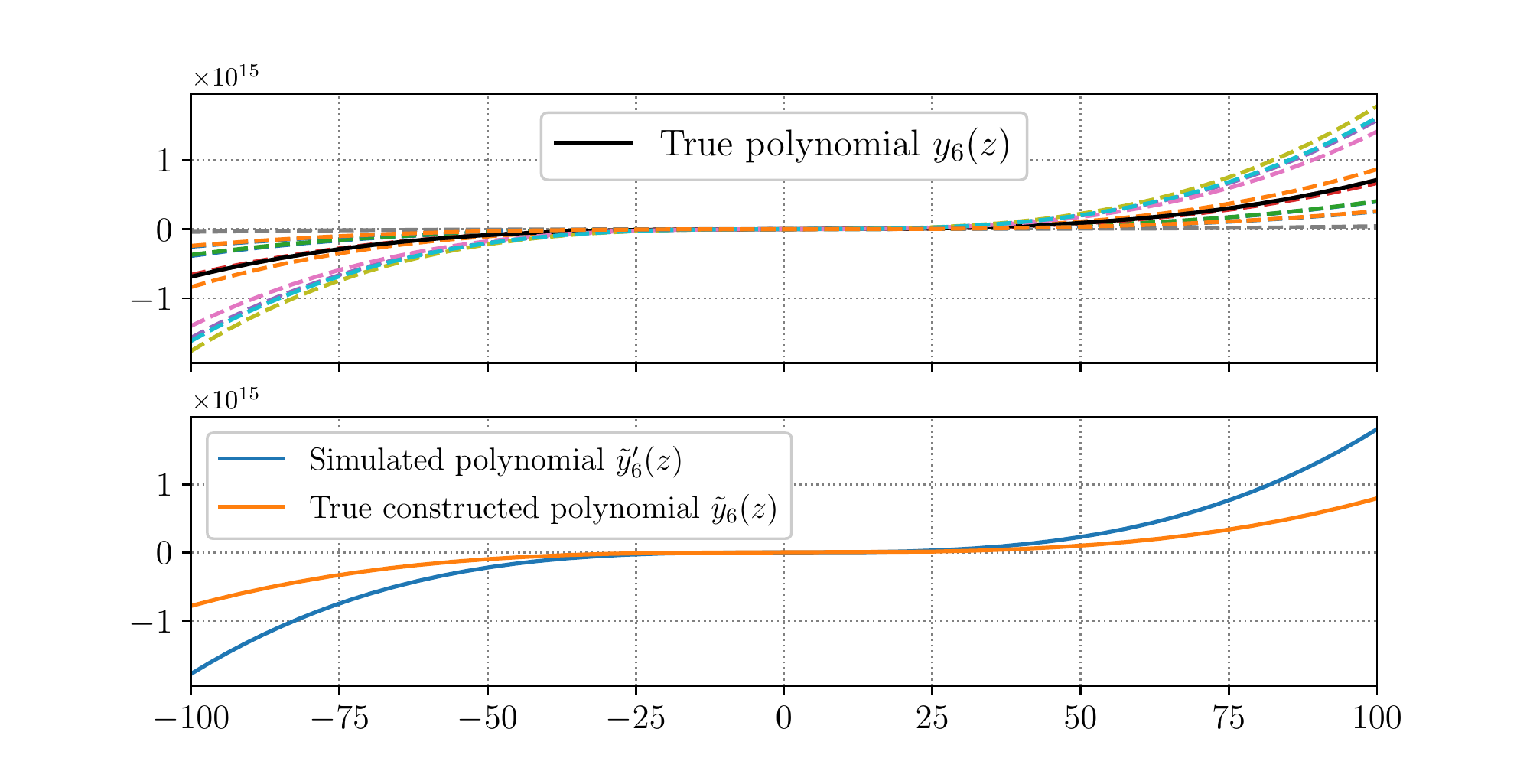} 
\vspace{-3mm}
 \caption{Polynomials simulated by a simulator to achieve computational indistinguishability among agents}
    \label{fig8}
    \vspace{-2mm}\end{figure} shows the existence of a simulator that can generate true polynomial $y_6(z)$ and simulated polynomials $y_i^{\prime}(z)$ (dashed lines), $\forall i=1,\ldots,n, i \neq 6$, such that $(\bm{\pi}_6 \tilde{\bm{P}})^{\prime}= \bm{\pi}_6 \tilde{\bm{P}}$. Therefore, the computational indistinguishability $\tilde{y}_6^{\prime} (\alpha_j) \stackrel{c}{\equiv} \tilde{y}_6 (\alpha_j) ,\forall j=1,\ldots,c$ is satisfied at any iteration and any time slot, and herein $\bm{\pi}_1 \tilde{\bm{P}},\ldots,\bm{\pi}_n \tilde{\bm{P}}$  can be securely computed among buses and the $i$th bus can only know the information contained in its own view $\operatorname{VIEW}_{i}$.



\section{Conclusion}
\label{Conclusion}
This paper proposed a novel decentralized privacy-preserving algorithm with cloud computing architecture for DER control in distribution networks. The DER control problem was formulated into a constrained optimization problem with the objectives of minimizing the line loss, PV curtailment cost, and ESS degradation cost. By integrating SS into the decentralized PGM, the proposed approach achieved privacy preservation for DER owners' private data, including the DERs' generation, consumption and daily electricity usage. The security of the proposed approach was proved rigorously with privacy guarantees and analyses against honest-but-curious agents and external eavesdroppers. Simulation results verified the applicability of the proposed approach on the modified IEEE 13-bus test feeder with controllable ESSs and solar PVs. Moreover, the designed methodology can be readily used in general large-scale decentralized optimization problems in the context of privacy preservation provisions.




 \begin{appendices}
\section{Derivation of the PGM Updates} \label{Derivation_proof} 

 We take the IEEE 13-bus test feeder in Fig. \ref{13bus} for example to illustrate the derivation of subgradients in \eqref{22}. 
To prove \eqref{22a}, we firstly consider the subgradient of the power loss minimization objective, the active power loss is
\begin{align}
    f_1(\bm{p}_1^{g},\ldots,\bm{p}_n^{g}) &= \delta_1 \sum_{l_{ij} \in \mathbb{L}} r_{ij} \left
(\frac{\| \bm{\mathcal{P}}_{ij} \|^2_2}{V_0^2}\right) =  \frac{\delta_1\bar{r}}{V_0^2} \sum_{\iota \in \mathbb{L}} \| \bm{\mathcal{P}}_{\iota} \|^2_2 \nonumber\\
&= \frac{\bar{\delta}_1}{2} \sum_{\iota \in \mathbb{L}} \| \bm{\mathcal{P}}_{\iota} \|^2_2.
\label{36sss}
\end{align}
Take \eqref{global_constraint} into \eqref{36sss}, we have
\begin{equation}
    f_1(\bm{p}_1^{g},\ldots,\bm{p}_n^{g}) =  \frac{\bar{\delta}_1}{2} \sum_{\iota \in \mathbb{L}} \| \Tilde{\bm{Z}}_{\iota}\Tilde{\bm{P}} \|^2_2. \label{54sss}
\end{equation}
Without loss of generality, assume the $\nu$th PV with decision variable  $\tilde{\bm{p}}_{\nu}$ is connected at bus $i$, we have 
\begin{align}
    \nabla_{\tilde{\bm{p}}_{\nu}} \mathcal{L}(\cdot) &= \delta_1 \nabla_{\tilde{\bm{p}}_{\nu}} f_1(\bm{p}_1^{g},\ldots,\bm{p}_n^{g}) +\delta_2 \nabla_{\tilde{\bm{p}}_{\nu}} f_2( \tilde{\bm{p}}_{\nu}) \nonumber\\
    &{+} \sum_{\iota=1}^{L} \nabla_{\tilde{\bm{p}}_{\nu}} \bm{\mu}_{u\iota}^{\mathsf{T}} (\Tilde{\bm{Z}}_{\iota}\Tilde{\bm{P}}{-}\overline{\bm{\mathcal{P}}}_{\iota})
    {-} \sum_{\iota=1}^{L} \nabla_{\tilde{\bm{p}}_{\nu}} \bm{\mu}_{l\iota}^{\mathsf{T}} \Tilde{\bm{Z}}_{\iota}\Tilde{\bm{P}}. \label{38s}
\end{align}
Substitute \eqref{16s} and \eqref{36sss} into the first term of \eqref{38s}, we have 
\begin{align}
   \delta_1 \nabla_{\tilde{\bm{p}}_{\nu}} f_1(\cdot) &= 
 \frac{\bar{\delta}_1}{2} \nabla_{\tilde{\bm{p}}_{\nu}}  \sum_{\iota \in \mathbb{L}} \| \Tilde{\bm{Z}}_{\iota}\Tilde{\bm{P}} \|^2_2 \nonumber\\
    &= \bar{\delta}_1 \sum_{\iota \in \mathbb{L}} \left( \nabla_{\tilde{\bm{p}}_{\nu}} \Tilde{\bm{Z}}_{\iota}  \sum_{\hat{\imath}=1}^{n}\bm{\Delta}_{\hat{\imath}} \tilde{\bm{p}}_{\hat{\imath}} \right)  \left( \Tilde{\bm{Z}}_{\iota}\Tilde{\bm{P}}\right) \nonumber\\
    &= \bar{\delta}_1 \sum_{\iota \in \mathbb{L}} \left(  \Tilde{\bm{Z}}_{\iota}  \bm{\Delta}_i \right)^{\mathsf{T}} \left( \Tilde{\bm{Z}}_{\iota}\Tilde{\bm{P}}\right). \label{39s}
\end{align}
Take the subgradient of \eqref{6ss},  the second term in \eqref{38s} becomes
\begin{equation}
    \delta_2\nabla_{\tilde{\bm{p}}_{\nu}} f_2(\tilde{\bm{p}}_{\nu}) = \delta_2\nabla_{\tilde{\bm{p}}_{\nu}} \| \tilde{\bm{p}}_{\nu} - \bm{\overline{p}}_{\nu}^{v}\|^2_2 = 2 \delta_2 \left(\tilde{\bm{p}}_{\nu} - \bm{\overline{p}}_{\nu}^{v}\right). \label{40s}
\end{equation}

Then, substitute \eqref{16s} into the third term of \eqref{38s} on the right hand side, we  have 
\begin{align}
\sum_{\iota=1}^{L} \nabla_{\tilde{\bm{p}}_{\nu}} \bm{\mu}_{u\iota}^{\mathsf{T}} (\Tilde{\bm{Z}}_{\iota}\Tilde{\bm{P}}-\overline{\bm{\mathcal{P}}}_{\iota}) 
&=  \sum_{\iota=1}^{L} \nabla_{\tilde{\bm{p}}_{\nu}} \bm{\mu}_{u\iota}^{\mathsf{T}} \Tilde{\bm{Z}}_{\iota}(\sum_{\hat{\imath}=1}^{n}\bm{\Delta}_{\hat{\imath}} \tilde{\bm{p}}_{\hat{\imath}} )\nonumber\\
&=  \sum_{\iota=1}^{L}  {(\Tilde{\bm{Z}}_{\iota}\bm{\Delta}_i)}^{\mathsf{T}} \bm{\mu}_{u\iota}. \label{41s}
\end{align}
Similarly, the last term of \eqref{38s} can be readily obtained as \vspace{-2mm}
\begin{equation}
-\sum_{\iota=1}^{L} \nabla_{\tilde{\bm{p}}_{\nu}} \bm{\mu}_{l\iota}^{\mathsf{T}} (\Tilde{\bm{Z}}_{\iota}\Tilde{\bm{P}}) = - \sum_{\iota=1}^{L}  {(\Tilde{\bm{Z}}_{\iota}\bm{\Delta}_i)}^{\mathsf{T}} \bm{\mu}_{l\iota}. \label{42s}
\end{equation}

Finally, by substituting \eqref{39s}, \eqref{40s}, \eqref{41s}, \eqref{42s} into \eqref{38s}, \eqref{22a} is readily proved. Following similar lines, subgradients of the primal variable $\hat{\bm{p}}_{\sigma}$ in \eqref{22b} can be readily proved. 

\section{Proof of \textbf{Theorem 2}} \label{Theorem2_proof}


\emph{Proof}: To prove the correctness of \textbf{Algorithm \ref{alg_1}}, we show that the proposed method has the same primal and dual solutions as the non-privacy PGM. 
Recall that the $u$th cloud multiplies the received $n$ outputs by the elements of $\bm{\pi}_i$ according to \eqref{36ss}, it yields
\begin{equation} \label{51}
\left\{\begin{array}{l}\bm{\pi}_i(1) y_{1}(\alpha_u) = \bm{\pi}_i(1) \left(\omega_{1}+a_{1,1} \alpha_{u} + \cdots + a_{1,d-1} \alpha_{u}^{d-1} \right)   \\ 
\vdots
\\
\bm{\pi}_i(n) y_{n}(\alpha_u) = \bm{\pi}_i(n) \left( \omega_{n} + a_{n,1} \alpha_{u} + \cdots + a_{n,d-1} \alpha_{u}^{d-1} \right) \\
\end{array}\right.
\end{equation}
Then, the aggregated outputs $\sum_{\hat{\imath}=1}^{n} \bm{\pi}_i(\hat{\imath}) y_{\hat{\imath}}(\alpha_u)$ in \eqref{36ss} can be obtained by summing the left hand side of \eqref{51}. Therefore, in total $c$ pairs of shares from all clouds as in \eqref{38ssss}  can be seen as the inputs and outputs of a polynomial 
\begin{equation} \label{53}
\tilde{y}(z) = \sum_{\hat{\imath}=1}^{n}
\bm{\pi}_i(\hat{\imath}) \omega_{\hat{\imath}}
+\tilde{a}_{1}z + \cdots + \tilde{a}_{d-1} z^{d-1}
\end{equation}
where $\tilde{a}_{\hat{\jmath}} = \sum_{\hat{\imath}=1}^{n}
\bm{\pi}_i(\hat{\imath}) a_{\hat{\imath},\hat{\jmath}}, \hat{\jmath}=1,\ldots,d-1$ and $\sum_{\hat{\imath}=1}^{n}
\bm{\pi}_i(\hat{\imath}) \omega_{\hat{\imath}}$ is exactly  $\bm{\pi}_{i}\tilde{\bm{P}}$. Then, the aggregated secret $\bm{\pi}_{i}\tilde{\bm{P}}$ can be readily retrieved by using $c$ pairs of shares in \eqref{37ss} since $d \leq c$, as stated by \textbf{Theorem 1}.

\section{Proof of \textbf{Proposition 1}} \label{Proposition_1_proof}

\textit{Proof}:  Under the collusion of $d-1$ clouds, they can construct the following set of equations
\begin{equation} \label{39sss}
\left\{\begin{array}{l}\tilde{y}_i(\alpha_1) = \tilde{\omega}
+\tilde{a}_{i,1} \alpha_1 + \cdots + \tilde{a}_{i,d-1} \alpha_1^{d-1}\\ 
\vdots
\\\tilde{y}_i(\alpha_{d-1}) = \tilde{\omega} +\tilde{a}_{i,1} \alpha_{d-1} + \cdots + \tilde{a}_{i,d-1} \alpha_{d-1}^{d-1} \\
\end{array}\right.
\end{equation}
where $\tilde{y}_i(z)$ is defined in \eqref{43} and $\tilde{\omega} = \bm{\pi}_{i}\tilde{\bm{P}}$. In \eqref{39sss}, $\tilde{a}_{i,\imath}$, $\forall \imath= 1,\ldots,d-1$ and $\Tilde{\omega}$ are unknown, therefore the $d-1$ clouds can yield in total $d-1$ equations yet $d$ unknowns that leads to underdetermined solutions.

\section{Proof of \textbf{Theorem 3}} \label{Theorem_3_proof}

 \emph{Proof}: To prove the privacy preservation of \textbf{Algorithm \ref{alg_1}} against honest-but-curious agents, we aim at verifying that whatever an honest-but-curious agent receives can be efficiently simulated. That being said, the honest-but-curious agent cannot retrieve useful information from  others using the received data because it cannot distinguish the received data from its own. During the $\ell$th iteration of executing \textbf{Algorithm \ref{alg_1}}, the view of  bus $i$ can be described via
\begin{align}
   \operatorname{VIEW}_{i} &= \{\alpha_1, \ldots,  \alpha_c, \theta, \gamma, \zeta, \bm{\pi}_i \tilde{\bm{P}},  y_i(z), \omega_i, \Bar{\mathcal{A}}_{i}, \nonumber \\ 
   &\qquad \quad \tilde{y}_i(\alpha_j), \forall j=1,\ldots,c, \mathcal{C}_p, \mathcal{C}_d \}.
\end{align}

Based on \textbf{Definition 2}, we need to prove the existence of a polynomial-time algorithm, denoted as simulator $\mathcal{S}$, that can simulate $\operatorname{VIEW}_i$ using the data of agent $i$, i.e., 
\begin{equation}
\mathcal{S}({\Xi}_i) \stackrel{c}{\equiv} \operatorname{VIEW}_i
\label{42ss}
\end{equation}
where $\Xi_i \triangleq \{\alpha_1, \ldots,  \alpha_c, \theta, \gamma, \zeta, \bm{\pi}_i \tilde{\bm{P}},  y_i(z), \omega_i, \Bar{\mathcal{A}}_{i}, \tilde{y}_i(\alpha_j)$, $\forall j=1,\ldots,c, \mathcal{C}_p, \mathcal{C}_d \} $ denotes the set of data that agent $i$ has access to. Manifesting  \eqref{42ss} indicates that whatever agent $i$ receives can be efficiently reconstructed based on its own knowledge  $\Xi_i$. To this end, the simulator is required to generate $\tilde{y}_i^{\prime} (\alpha_j)$,$\forall j=1,\ldots,c$ that satisfy 
\begin{equation}
\tilde{y}_i^{\prime} (\alpha_j) \stackrel{c}{\equiv} \tilde{y}_i (\alpha_j) ,\forall j=1,\ldots,c.
\label{55ssp}
\end{equation}

To achieve this goal, the simulator firstly generates secrets $w^{\prime}_{j \neq i} \in \mathbb{E}$ of other agents such that 
\begin{equation}
   \bm{\pi}_i \tilde{\bm{P}} = w_i + \sum
   _{j\neq i} w_j^{\prime}.
\end{equation}
Then it generates a set of random polynomials as in \eqref{34s} to obtain $y_j^{\prime}(z), \forall j \neq i$ with $w_j^{\prime}, \forall j \neq i$ as the corresponding constant terms, i.e.,
\begin{subnumcases}{\label{59sss}}
\begin{aligned}
    &y_i(z) = w_i + a_{i,1}z + \cdots + a_{i,d-1}z^{d-1}
\end{aligned}\label{59sss:a}\\
\begin{aligned}
& y_j^{\prime}(z) = w_j^{\prime} + a_{i,1}^{\prime}z + \cdots + a_{i,d-1}^{\prime}z^{d-1}, \forall j \neq i.
\end{aligned} \label{59sss:b}
\end{subnumcases}

Consequently, the simulator can use $\{\alpha_1, \ldots, \alpha_c\}$ as inputs for  \eqref{59sss} and obtain 
\begin{equation}
    \Tilde{\mathcal{A}}_{i}^{\prime} = \left\{
    \alpha_{\hat{\jmath}}, y_{i}(\alpha_{\hat{\jmath}})+\sum_{j \neq i} y_{j}^{\prime}(\alpha_{\hat{\jmath}}),\forall,{\hat{\jmath}}=1,\ldots,c \right\}.
    \label{37sse}
\end{equation}

By \textbf{Theorem 1} and \textbf{Theorem 2}, the shares in \eqref{37sse} can be used to construct a new polynomial in the form of
\begin{equation}
\tilde{y}_i^{\prime}(x) = (\bm{\pi}_i \tilde{\bm{P}})^{\prime} + \tilde{a}_{i,1}^{\prime}z + \cdots + \tilde{a}_{i,d-1}^{\prime}z^{d-1} \label{45ss}
\end{equation}
where $(\bm{\pi}_i \tilde{\bm{P}})^{\prime} = \bm{\pi}_i \tilde{\bm{P}}$. Therefore, \eqref{55ssp} and \eqref{42ss} hold, by \textbf{Definition 2}, \textbf{Algorithm 1} securely computes $\bm{\pi}_1 \tilde{\bm{P}},\ldots,\bm{\pi}_n \tilde{\bm{P}}$ between the agents.

In what follows, we prove the privacy preservation of \textbf{Algorithm \ref{alg_1}} against external eavesdroppers. Under \textbf{Assumption 1}, assume agent 1 is safe from external eavesdroppers, by wiretapping any other agents' communication channels, an external eavesdropper can at most have access to 
\begin{equation}
   \Xi_e {=} \left\{\alpha_1{,} {\ldots}, \alpha_c{,} y_i(\alpha_u), \Bar{\mathcal{A}}_{u,i}, \forall i {=} 2{,}{\ldots},n{,} u {=} 1,{\ldots},c \right\}. \label{60st}
\end{equation}
Since \eqref{60st} is insufficient  to formulate \eqref{37ss}, the external eavesdropper is incapable of inferring either $y_i(z)$'s or $\tilde{y}_i^{\prime}(z)$'s, i.e., unable to infer agents' private information $\bm{p}_i$'s or the aggregated message $\bm{\pi}_i \tilde{\bm{P}}$'s.

\end{appendices}

\bibliographystyle{IEEEtran}

\bibliography{bibliography}

\begin{thebibliography}{10}
\providecommand{\url}[1]{#1}
\csname url@samestyle\endcsname
\providecommand{\newblock}{\relax}
\providecommand{\bibinfo}[2]{#2}
\providecommand{\BIBentrySTDinterwordspacing}{\spaceskip=0pt\relax}
\providecommand{\BIBentryALTinterwordstretchfactor}{4}
\providecommand{\BIBentryALTinterwordspacing}{\spaceskip=\fontdimen2\font plus
\BIBentryALTinterwordstretchfactor\fontdimen3\font minus
  \fontdimen4\font\relax}
\providecommand{\BIBforeignlanguage}[2]{{%
\expandafter\ifx\csname l@#1\endcsname\relax
\typeout{** WARNING: IEEEtran.bst: No hyphenation pattern has been}%
\typeout{** loaded for the language `#1'. Using the pattern for}%
\typeout{** the default language instead.}%
\else
\language=\csname l@#1\endcsname
\fi
#2}}
\providecommand{\BIBdecl}{\relax}
\BIBdecl

\bibitem{derancillary}
J.~Campbell, ``Ancillary services provided from {DER},'' Oak Ridge National
  Lab, Oak Ridge, TN, United States, Tech. Rep., 2005.

\bibitem{abb_aguero2017modernizing}
J.~R. Aguero, E.~Takayesu, D.~Novosel, and R.~Masiello, ``Modernizing the grid:
  {C}hallenges and opportunities for a sustainable future,'' \emph{IEEE Power
  Energy Mag.}, vol.~15, no.~3, pp. 74--83, 2017.

\bibitem{abb_song2012operation}
I.-K. Song, W.-W. Jung, J.-Y. Kim, S.-Y. Yun, J.-H. Choi, and S.-J. Ahn,
  ``Operation schemes of smart distribution networks with distributed energy
  resources for loss reduction and service restoration,'' \emph{IEEE Trans.
  Smart Grid}, vol.~4, no.~1, pp. 367--374, 2012.

\bibitem{abb_molzahn2017survey}
D.~K. Molzahn, F.~D{\"o}rfler, H.~Sandberg, S.~H. Low, S.~Chakrabarti,
  R.~Baldick, and J.~Lavaei, ``A survey of distributed optimization and control
  algorithms for electric power systems,'' \emph{IEEE Trans. Smart Grid},
  vol.~8, no.~6, pp. 2941--2962, 2017.

\bibitem{abb_zeraati2016distributed}
M.~Zeraati, M.~E.~H. Golshan, and J.~M. Guerrero, ``Distributed control of
  battery energy storage systems for voltage regulation in distribution
  networks with high {PV} penetration,'' \emph{IEEE Trans. Smart Grid}, vol.~9,
  no.~4, pp. 3582--3593, 2016.

\bibitem{abb_zhang2021distributed}
Q.~Zhang, Y.~Guo, Z.~Wang, and F.~Bu, ``Distributed optimal conservation
  voltage reduction in integrated primary-secondary distribution systems,''
  \emph{IEEE Trans. Smart Grid}, vol.~12, no.~5, pp. 3889--3900, 2021.

\bibitem{abb_pan2021distributed}
Y.~Pan, A.~Sangwongwanich, Y.~Yang, and F.~Blaabjerg, ``Distributed control of
  islanded series {PV}-battery-hybrid systems with low communication burden,''
  \emph{IEEE Trans. Power Electron.}, vol.~36, no.~9, pp. 10\,199--10\,213,
  2021.

\bibitem{abb_navidi2018two}
T.~Navidi, A.~El~Gamal, and R.~Rajagopal, ``A two-layer decentralized control
  architecture for {DER} coordination,'' in \emph{Proc. IEEE Conf. Decis.
  Control}, Miami, FL, USA, Dec. 17-29 2018, pp. 6019--6024.

\bibitem{abb_lin2017decentralized}
W.~Lin and E.~Bitar, ``Decentralized stochastic control of distributed energy
  resources,'' \emph{IEEE Trans. Power Syst.}, vol.~33, no.~1, pp. 888--900,
  2017.

\bibitem{abb_huo2022two}
X.~Huo and M.~Liu, ``Two-facet scalable cooperative optimization of multi-agent
  systems in the networked environment,'' \emph{IEEE Trans. Control Syst.
  Technol.}, vol.~30, no.~6, pp. 2317--2332, 2022.

\bibitem{abb_dwork2006calibrating}
C.~Dwork, F.~McSherry, K.~Nissim, and A.~Smith, ``Calibrating noise to
  sensitivity in private data analysis,'' in \emph{Proc. Theory Cryptogr.
  Conf.}, New York, NY, USA, Mar. 4-7 2006, pp. 265--284.

\bibitem{abb_dong2018privacy}
J.~Dong, T.~Kuruganti, S.~Djouadi, M.~Olama, and Y.~Xue, ``Privacy-preserving
  aggregation of controllable loads to compensate fluctuations in solar
  power,'' in \emph{Proc. IEEE Electron. Power Grid}, Charleston, SC, USA, Nov.
  12-14 2018, pp. 1--5.

\bibitem{abb_han2016differentially}
S.~Han, U.~Topcu, and G.~J. Pappas, ``Differentially private distributed
  constrained optimization,'' \emph{IEEE Trans. Autom. Control}, vol.~62,
  no.~1, pp. 50--64, 2016.

\bibitem{abb_gough2021preserving}
M.~B. Gough, S.~F. Santos, T.~AlSkaif, M.~S. Javadi, R.~Castro, and J.~P.
  Catal{\~a}o, ``Preserving privacy of smart meter data in a smart grid
  environment,'' \emph{IEEE Trans. Ind. Inform.}, vol.~18, no.~1, pp. 707--718,
  2021.

\bibitem{abb_lu2012eppa}
R.~Lu, X.~Liang, X.~Li, X.~Lin, and X.~Shen, ``{EPPA}: An efficient and
  privacy-preserving aggregation scheme for secure smart grid communications,''
  \emph{IEEE Trans. Parallel Distrib. Syst.}, vol.~23, no.~9, pp. 1621--1631,
  2012.

\bibitem{abb_mohammadali2021privacy}
A.~Mohammadali and M.~S. Haghighi, ``A privacy-preserving homomorphic scheme
  with multiple dimensions and fault tolerance for metering data aggregation in
  smart grid,'' \emph{IEEE Trans. Smart Grid}, vol.~12, no.~6, pp. 5212--5220,
  2021.

\bibitem{abb_cheng2021homomorphic}
Z.~Cheng, F.~Ye, X.~Cao, and M.-Y. Chow, ``A homomorphic encryption-based
  private collaborative distributed energy management system,'' \emph{IEEE
  Trans. Smart Grid}, vol.~12, no.~6, pp. 5233--5243, 2021.

\bibitem{abb_wang2018privacy}
S.~Wang, Q.~Hu, Y.~Sun, and J.~Huang, ``Privacy preservation in location-based
  services,'' \emph{IEEE Commun. Mag.}, vol.~56, no.~3, pp. 134--140, 2018.

\bibitem{abb_gilad2019secure}
R.~Gilad-Bachrach, K.~Laine, K.~Lauter, P.~Rindal, and M.~Rosulek, ``Secure
  data exchange: {A} marketplace in the cloud,'' in \emph{Proc. ACM SIGSAC
  Conf. Cloud Comput. Security Workshop}, London, UK, Nov. 11 2019, pp.
  117--128.

\bibitem{abb_shamir1979share}
A.~Shamir, ``How to share a secret,'' \emph{Commun. {ACM}}, vol.~22, no.~11,
  pp. 612--613, 1979.

\bibitem{nabil2019ppetd}
M.~Nabil, M.~Ismail, M.~M. Mahmoud, W.~Alasmary, and E.~Serpedin, ``{PPETD}:
  Privacy-preserving electricity theft detection scheme with load monitoring
  and billing for {AMI} networks,'' \emph{IEEE Access}, vol.~7, pp.
  96\,334--96\,348, 2019.

\bibitem{abb_huo2022distributed}
X.~Huo and M.~Liu, ``Distributed privacy-preserving electric vehicle charging
  control based on secret sharing,'' \emph{Electr. Power Syst. Res.}, vol. 211,
  p. 108357, 2022.

\bibitem{abb_baran1989optimal}
M.~Baran and F.~F. Wu, ``Optimal sizing of capacitors placed on a radial
  distribution system,'' \emph{IEEE Trans. Power Deliv.}, vol.~4, no.~1, pp.
  735--743, 1989.

\bibitem{abb_farivar2013equilibrium}
M.~Farivar, L.~Chen, and S.~Low, ``Equilibrium and dynamics of local voltage
  control in distribution systems,'' in \emph{Proc. IEEE Conf. Decis. Control},
  Florence, Italy, Dec. 10-13 2013, pp. 4329--4334.

\bibitem{abb_li2019distributed}
J.~Li, Z.~Xu, J.~Zhao, and C.~Zhang, ``Distributed online voltage control in
  active distribution networks considering {PV} curtailment,'' \emph{IEEE
  Trans. Ind. Inform.}, vol.~15, no.~10, pp. 5519--5530, 2019.

\bibitem{abb_forman2013optimization}
J.~Forman, J.~Stein, and H.~Fathy, ``Optimization of dynamic battery parameter
  characterization experiments via differential evolution,'' in \emph{Proc. Am.
  Control Conf.}, Washington, DC, USA, Jun. 17-19 2013, pp. 867--874.

\bibitem{abb_daru2019encrypted}
M.~S. Daru and T.~Jager, ``Encrypted cloud-based control using secret sharing
  with one-time pads,'' in \emph{Proc. IEEE Conf. Decis. Control}, Nice,
  France, Dec. 11-13 2019, pp. 7215--7221.

\bibitem{abb_humpherys2020foundations}
J.~Humpherys and T.~J. Jarvis, \emph{Foundations of Applied Mathematics, Volume
  I: Mathematical Analysis}.\hskip 1em plus 0.5em minus 0.4em\relax Soc. Ind.
  Appl. Math, 2020.

\bibitem{abb_goldreich2009foundations}
O.~Goldreich, \emph{Foundations of Cryptography: Volume 2, Basic
  Applications}.\hskip 1em plus 0.5em minus 0.4em\relax Cambridge Univ. Press,
  2009.

\bibitem{abb_evans2018pragmatic}
D.~Evans, V.~Kolesnikov, and M.~Rosulek, ``A pragmatic introduction to secure
  multi-party computation,'' \emph{Found. Trends Privacy Security}, vol.~2, no.
  2-3, pp. 70--246, 2018.

\bibitem{goldreich1998secure}
O.~Goldreich, ``Secure multi-party computation,'' \emph{Manuscript. Preliminary
  Version}, vol.~78, p. 110, 1998.

\bibitem{abb_liu2017decentralized}
M.~Liu, P.~K. Phanivong, Y.~Shi, and D.~S. Callaway, ``Decentralized charging
  control of electric vehicles in residential distribution networks,''
  \emph{IEEE Trans. Control Syst. Technol.}, vol.~27, no.~1, pp. 266--281,
  2019.

\bibitem{ESS}
\BIBentryALTinterwordspacing
{National Renewable Energy Laboratory}. Residential battery storage. [Online].
  Available:
  \url{https://atb.nrel.gov/electricity/2021/residential_battery_storage}
\BIBentrySTDinterwordspacing

\bibitem{CISO}
\BIBentryALTinterwordspacing
{U.S. Energy Information Administration}. Electric power annual. [Online].
  Available: \url{https://www.eia.gov/todayinenergy/detail.php?id=49276}
\BIBentrySTDinterwordspacing

\end{thebibliography}


\end{document}